\documentclass[final,onefignum,onetabnum]{siamltex1213}

\usepackage{cite}
\usepackage{amsfonts,amssymb,amsmath}
\usepackage{graphicx}

\newtheorem{alg}{Algorithm}
\newcommand{\algo}[1]{\medskip%
  \noindent \parbox{1.\textwidth}{\begin{alg}%
     \rm\ \\*[1ex] \fbox{\ \parbox{0.95\textwidth}{ \rule[1ex]{0pt}{2ex}  #1 \rule[0ex]{0pt}{2ex}  }\ }%
     \end{alg}}\bigskip}
\newcommand{\cT}{\mathcal{T}}
\newcommand{\ZZ}{\mathbb{Z}}
\newcommand{\RR}{\mathbb{R}}
\newcommand{\I}{\mathcal{I}}
\newcommand{\J}{\mathcal{J}}
\newcommand{\spn}{\mathrm{span}}
\newcommand{\NN}{\mathbb{N}}
\newcommand{\No}{\mathbb{N}_0}

\title{Fast Algorithms for Adaptive Free-Knot Spline Approximation Using Non-Uniform Biorthogonal Spline Wavelets
\thanks{This work has been 
partly supported by the ENIAC research project ARTEMOS under grant 829397 and
the FWF under grant P22549.}} 

\author{Kai Bittner\footnotemark[2]\ \footnotemark[3]
\and Hans Georg Brachtendorf\footnotemark[2]\ \footnotemark[4] }

\begin{document}
\maketitle
\renewcommand{\thefootnote}{\fnsymbol{footnote}}

\footnotetext[2]{University of Applied Sciences Upper Austria, Softwarepark 11, A-4232 Hagenberg, Austria}
\footnotetext[3]{Kai.Bittner@fh-hagenberg.at}
\footnotetext[4]{Hans-Georg.Brachtendorf@fh-hagenberg.at}

\renewcommand{\thefootnote}{\arabic{footnote}}

\slugger{sisc}{xxxx}{xx}{x}{x--x}

\begin{abstract}
New algorithms for fast wavelet transforms with biorthogonal spline wavelets on nonuniform grids are
presented. In contrary to classical wavelet transforms, the algorithms are not based on filter coefficients,
but on algorithms for B-spline expansions (differentiation, Oslo algorithm, etc.). Due to inherent properties
of the spline wavelets, the algorithm can be modified for spline grid refinement or coarsening. The performance 
of the algorithms is demonstrated by numerical tests of the adaptive spline methods in circuit simulation.
\end{abstract}

\begin{keywords}
Splines, spline wavelets, free knot spline approximation 
\end{keywords}

\begin{AMS}
65D07, 41A15, 65T60, 42C40
\end{AMS}

\pagestyle{myheadings}
\thispagestyle{plain}
\markboth{K.~BITTNER AND H.G.~BRACHTENDORF}{FAST ALGORITHMS FOR NON-UNIFORM SPLINE WAVELETS}

\section{Introduction}

Since the dawn of wavelet theory spline wavelets have been always of particular interest.
This includes orthogonal spline wavelets \cite{Battle87,Lemar88}, semi-orthogonal spline wavelets 
\cite{ChWa92a,ChWa92b} as well as
biorthogonal wavelets \cite{CDF}. An exceptional property of spline wavelets is that they possess an explicit representation
(in terms of piecewise polynomials), while most other wavelets of interest are only described by their two-scale relation.
This permits extra flexibility, e.g.\ for the construction of wavelets on the interval \cite{ChQu92,DKU2}
or the evaluation of non-linear mappings of wavelet expansions \cite{DSX00,CDD03,BiUr04}.

The wavelet constructions above are based on equidistant spline knots. However, since a spline can 
be defined for any given grid, it arises the question if spline wavelet constructions on nonuniform grids are possible.
There have been several publications on semi-orthogonal spline wavelets on non-uniform grids 
(see e.g.~\cite{BuhMic92,KazElb97,LyMoQu01}). Although these wavelets do not have a sparse decomposition relation,
there are fast algorithms, which solve a banded linear system involving the reconstruction coefficients (cf.~\cite{QW94}).
However, in some cases finite decomposition relations as they appear for biorthogonal spline wavelets \cite{CDF,DKU2}
may be of interest. A constructive proof of the existence of such spline wavelets in the nonuniform setting was given
in \cite{DahMic93}, based on results of banded matrices with banded inverses. However, no algorithms based on this
approach are provided.

In \cite{Bit05b} we have given sufficient and necessary conditions for the existence of finite reconstruction
and decomposition relations. However, an algorithm based on this relation requires the computation of many coefficients,
which is time and memory consuming. Thus, we have developed a direct approach based on known properties of spline
functions. Here we will present
algorithms for the fast wavelet transform for non-uniform spline wavelets, which represent a generalization
of the biorthogonal spline wavelets from \cite{CDF}. By a small modification we obtain also an algorithm for 
adaptive knot removal, which permits to reduce the size of a spline representation with the approximation error 
under control. Furthermore, we introduce an adaptive spline approximation method with adaptive grid
refinement.
These algorithms were used to develop a wavelet based adaptive method for circuit simulation
\cite{BiDau10a,BiDau10b,BiBra12b,BiBra14b}. In this context the use of wavelets on non-uniform grids 
appears to be much more suitable than a method based on uniform wavelets.
Here we will give a complete description of the wavelet transforms.

In \S\ref{splines} we will recall basic properties of spline functions and B-splines, which are important for
our approach. The wavelet bases for non-uniform grids will be introduced in \S\ref{wavelets}.
The fast decomposition and reconstruction algorithm is described in \S\ref{algorithms}.
\S\ref{modify} discusses several modifications of our approach. In particular the restriction
to an interval and to periodic functions is of practical interest. Based on the fast decomposition algorithm,
methods for grid coarsening and refinement are introduced in \S\ref{wave_adapt}. The performance
of this grid adaptation is demonstrated by a numerical test on a multirate circuit simulation problem
in \S\ref{num_test}.

\section{Spline functions and B-splines\label{splines}}

In this section we establish our notation and recall basic facts and algorithms for splines, which are
needed for our new algorithms. For a detailed introduction to splines we refer to \cite{Boor78,Sch81}.

For a given knot sequence $\cT=(t_i)_{i\in\ZZ}$ satisfying
\begin{equation}\label{knot_cond}
t_i< t_{i+1},
\end{equation}
the \emph{spline space} of order $m$ is defined as
$$
S_{m,\cT} =\Big\{ f\in C^{m-2}(\RR):~f\big|_{(t_i,t_{i+1}]}\in \Pi_{m-1}\Big\}.
$$
Here, $\Pi_{m-1}$ denotes the space of polynomials of degree less than $m$.
It is also possible to consider multiple knots, i.e., $t_i=\ldots=t_{i+\ell}$.
In this case the spline functions possess only $m-2-\ell$ continuous derivatives in $t_i$.
Since we will need multiple knots for spline wavelets on a closed interval we include
the corresponding modifications in the sequel.

A basis for $S_{m,\cT}$ is given by the B-splines $N_{m,k}=N^{\cT}_{m,k}$, 
$k\in\ZZ$, which are
the uniquely determined spline functions of minimal support $[t_k,t_{k+m}]$
satisfying $\sum_{k\in\ZZ} N_{m,k}(1)\equiv 1$.

Values of a spline function $s$ can be computed from the coefficients $c_k$ of its B-spline expansion
$s=\sum_{k\in\ZZ} c_k N_{m,k}$ by the de Boor algorithm. 
For $t\in (t_k,t_{k+1}]$ the value 
$s(t)=c^{(m-1)}_k$ is computed by the recursion
\begin{equation}
\label{deBoor}
c_k^{(\ell)}=\left\{\begin{array}{ll}
c_k,&\mbox{if }\ell=0,\\
\displaystyle
   \frac{(t-t_k)\, c_k^{(\ell-1)}+(t_{k+m-\ell}-t)\,c_{k-1}^{(\ell-1)}}{
              t_{k+m-\ell}-t_k},& \mbox{if }t\in (t_k,t_{k+m-\ell}],\\
0& \mbox{otherwise.}
\end{array}\right.
\end{equation}

Furthermore, derivatives of the B-splines are given by
\begin{equation}
\label{bsp_der}
N^\prime_{m,k}(t)=(m-1)\left(\frac{N_{m-1,k}(t)}{t_{k+m-1}-t_k}
               -\frac{N_{m-1,k+1}(t)}{t_{k+m}-t_{k+1}}\right).
\end{equation}
Therefore the derivative of the spline function 
$s=\sum_{k=0}^n c_k N_{m,k}$ is
given as
$s'=\sum_{k=0}^n c^\prime_k N_{m-1,k}$,
with
\begin{equation}
\label{spl_der}
c_k^\prime=\left\{\begin{array}{ll}
\displaystyle
   (m-1)\frac{c_k-c_{k-1}}{
              t_{k+m-1}-t_k}& \mbox{if }t_k<t_{k+m-1},\\
0& \mbox{otherwise}.
\end{array}\right.
\end{equation}

For our wavelet construction we have to deal with nested
spline spaces. Obviously, the spline space $S_{m,\Delta}$ is contained
in $S_{m,\cT}$ if and only if 
$\Delta=(\tau_k)_{k\in\ZZ}\subset\cT=(t_\ell)_{\ell\in\ZZ}$, i.e.,
$\cT$ is a \emph{refinement} of $\Delta$.

Any spline $s=\sum_{k=1}^N c^\Delta_k N_{m,k}^\Delta\in S_{m,\Delta}$ is also contained
in $S_{m,\cT}$ and has the B-spline representation
$s=\sum_{k=1}^n c^\cT_k N_{m,k}^\cT$ with respect to the finer grid.
The coefficients $c^\cT_k$ can be computed from the $c^\Delta_k$ by the Oslo algorithm \cite{CLR80, LyMo86} as follows.
For $t_i\in (\tau_j,\tau_{j+1}]$ the coefficient $c^\cT_i=c_{j,i}^{(m-1)}$ is determined by the recursion 
\begin{equation}\label{oslo}
c_{j,i}^{(k)}=\left\{\!\!\begin{array}{ll}
c^\Delta_j,&\mbox{if }k=0,\\
   \dfrac{(t_{i+m-k}-\tau_j) c_{j,i}^{(k-1)}+
       (\tau_{j+m-k}-t_{i+m-k})c_{j-1,i}^{(k-1)}}{\tau_{j+m-k}-\tau_j},\!\!& 
              \mbox{if\;}t_{i+m-k}\in(\tau_j,\tau_{j+m-k}],\\
0&\mbox{otherwise}.
\end{array}\right.\!\!
\end{equation}

The insertion of a single knot, i.e.,
$\cT=(\tau_1,\ldots,\tau_{j-1},t,\tau_j,\ldots,\tau_N)$ can also be done
by Boehm's knot insertion algorithm \cite{Boe80} 
\begin{equation}
\label{boehm}
c^\cT_i=\left\{\begin{array}{ll}
c^\Delta_i,&\mbox{if }i\le j-m\\
\displaystyle\frac{(\tau_{i+m}-t)\, c^\Delta_{i-1}
               +(t-\tau_i)\,c^\Delta_i}{\tau_{i+m}-\tau_i},&\mbox{if }j-m<i< j,\\
c^\Delta_{i-1},&\mbox{if }i\ge j.
\end{array}\right.
\end{equation}

\section{A class of spline wavelets for non-uniform grids\label{wavelets}}

Spline wavelets on nonuniform grids with finite decomposition and reconstruction sequences have been investigated in \cite{Bit05b}. To have an efficient an readable algorithm we consider a subclass of these spline wavelets,
which can be considered as  a generalization of spline wavelets
for uniform grids \cite{CDF}, as pointed out in \cite{Bit05a}. We will illustrate possible generalizations after establishing
our fast algorithm.

We start from a lattice of nested knot sets $\cT^j=(t_i^j)_{i\in\ZZ}$
satisfying $t^j_i<t^j_{i+1}$, where  either
$t_i^j=t_{2i}^{j+1}$ or
\begin{equation}\label{local_refin}
t_i^j=\left\{\begin{array}{ll}
t_{n^j_1+i}^{j+1},&\mbox{if }i<n^j_1,\\[1.7ex]
t_{2i}^{j+1},&\mbox{if }i=n^j_1,\ldots,n^j_2\\[1.2ex]
t_{n^j_2+i}^{j+1},&\mbox{if }i>n^j_2.
\end{array}\right.
\end{equation}
The classical setting would be $n^j_1 = -\infty$ and $n^j_2 = \infty$, while finite values
for $n^j_1$ and $n^j_2$ yield a local refinement of $\cT^j$, where
additional knots are only inserted in the interval $[t^0_{n^j_1},t^0_{n^j_2}]$.
Obviously, the corresponding spline spaces $V_j:=S_{m,\cT^j}$
are nested, i.e., $V_j \subset V_{j+1}$.  Therefore, we consider $\{V_j:~j\in\ZZ\}$ as generalized
Multiresolution Analysis with scaling functions 
$$
\varphi^j_k:=N^{\cT^j}_{m,k}.
$$
Then, we define wavelets as  
\begin{equation}\label{wave_bit}
\psi^j_k(t):=\alpha^j_k\;\frac{d^{\widetilde{m}}}{dt^{\widetilde{m}}} N^{\Xi^j_k}_{m+\widetilde{m},0}(t),
\quad k= n^j_1,\ldots,n^j_2-1.
\end{equation}
where
\begin{equation}\label{wave_knots}
\Xi^j_k=\big(\xi^j_{k,\ell}\big)_{\ell=0}^{m+\widetilde{m}}
:=\big(t^j_{k+1-\ell_1},t^j_{k+2-\ell_1},\ldots
   t^j_k,t^{j+1}_{2k+1},t^j_{k+1},\ldots,
   t^j_{k-1+\ell_2},t^j_{k+\ell_2}\big)
\end{equation}
with $\ell_1=\big\lfloor \frac{m+\widetilde{m}}{2}\big\rfloor$,
$\ell_2=\big\lceil \frac{m+\widetilde{m}}{2}\big\rceil$.
Note, that $t^{j+1}_{2k+1}$ is the only knot in $\Xi^j_k$ which is not
contained in $\cT^j$. Thus, each wavelet represents the contribution of one and only one knot from 
$\cT^{j+1}\setminus\cT^j$. This property will be a key ingredient for our fast wavelet transform algorithm.
The term $\alpha^j_k$ in (\ref{wave_bit}) denotes a normalization factor, which can e.g.\ be chosen according 
to a particular function space. In \S\ref{coarse} (Eq.~(\ref{normalization}) and
enclosing text) we give an example how this factor can be chosen in an application.

As demonstrated in \cite{Bit05b,Bit05a}, we have the following properties:
\begin{romannum}
\item $\psi^j_k\in V_{j+1}$ and $\psi^j_k\not\in V_j$, 
\item Compact support:\quad $\supp \psi^j_k = [t^j_{k+1-\ell_1}, t^j_{k+\ell_2}]$,
\item $\widetilde{m}$ Vanishing moments:\quad 
  $\int_\RR \psi^j_k(t)\,t^\ell\, dt = 0,\quad \ell=0,\ldots,\widetilde{m}-1$,
\item Finite reconstruction relations:
\begin{equation}\label{rekrel}
\varphi^j_k=\sum_{\ell=2k}^{2k+m} a^j_{k,\ell}\,\varphi^{j+1}_\ell\quad \mbox{and}\quad
\psi^j_k=\sum_{\ell=2(k+1-\ell_1)}^{2(k+\ell_2)-m} b^j_{k,\ell}\,\varphi^{j+1}_\ell,
\end{equation}
\item Finite decomposition relations:
\begin{equation}\label{dekrel}
\varphi^{j+1}_k=\sum\limits_{\ell\in\I^j_k} c^j_{k,\ell}\,\varphi^j_\ell
         +\sum_{\ell\in\J^j_k} d^j_{k,\ell}\,\psi^j_\ell\qquad \#\I^j_k<\widetilde{m}+\big\lceil\tfrac{m}{2}\big\rceil,
\quad\#\J^j_k<\big\lceil\tfrac{m}{2}\big\rceil,
\end{equation}
\item Generalization of existing wavelets: For uniform dyadic grids, i.e.~if $t^j_k=2^{-j}k$,  
and for $m+\tilde{m}$ even, the wavelets
coincide with the biorthogonal spline wavelets introduced in \cite{CDF}.
These wavelets are widely used in Numerical Analysis, since there exist well understood generalizations 
to intervals \cite{Bit05b,DKU2,Primbs10} and 
domains \cite{DS,KuSa05} and since spline wavelets are also optimal for the handling of non-linear 
problems \cite{BiUr04,BDS04,CDD03,DSX00}.
\end{romannum}
\bigskip

Most definitions of wavelet systems include also a stability condition (e.g.\ orthonormal basis, Riesz basis), which ensures 
that small disturbances of the coefficients cause only small deviations of the spline function and vice versa.
While proofs of Riesz stability exist in the uniform case, similar results are much harder to obtain for the nonuniform case.
While there is some numerical evidence that the wavelet transform are stable, the theoretical foundation will be the
target of future research. 

The relations (\ref{rekrel}) and (\ref{dekrel}) permit a fast wavelet transform.
The wavelet transforms (decomposition and reconstruction) are a change of basis between
$\{\varphi^{j+1}_k\}$ and $\{\varphi^j_k\}\cup\{\psi^j_k\}$, which can be done fast due to the finiteness of
the sums. To perform these transforms we need to know the coefficients $a^j_{k,\ell}$, $b^j_{k,\ell}$, 
$c^j_{k,\ell}$, and $d^j_{k,\ell}$. The computation of these coefficients is possible but costly in computation time
and storage. Therefore, we propose here a direct way based on the proofs in \cite{Bit05b,Bit05a}, which is more practicable.

\section{Algorithms for the fast wavelet transform\label{algorithms}}

\subsection{Fast wavelet decomposition\label{wv_deco}}

The wavelet decomposition algorithm, also called wavelet analysis or just
fast wavelet transform, computes from given spline coefficients
$c_k^{j+1}$ the coarse scale spline coefficients $c_k^j$ and the corresponding
wavelet coefficients $d_k^j$ satisfying 
\begin{equation}\label{two_scale}
s_{j+1}:=\sum_{k\in\ZZ} c^{j+1}_k \varphi^{j+1}_k = \sum_{k\in\ZZ} c^j_k \varphi^j_k + d^j_k \psi^j_k.
\end{equation}
In the sequel we assume without loss of generality that $j = 0$.
Furthermore, we assume that $s_1$ is compactly supported in $[t_0^0,t_n^0]$ such that all sums are finite,
which is reasonable from the computational point of view. Obviously, we can then assume that in (\ref{local_refin})
$n^0_1 = 0$ and $n^0_2 = n$.

\paragraph{Computation of the wavelet coefficients}
Any spline can be written in its truncated power representation 
\begin{equation}
\label{s_trunc1}
s_1(t) = \sum_{k\in\ZZ} a_k\; (t-t^1_k)_+^{m-1} ,
\end{equation}
where  $(t)_+^j=\chi_{(0,\infty)}(t)\, t^j$
are the truncated powers. Due to (\ref{wave_bit}) and (\ref{wave_knots}), the wavelets have the 
representation
\begin{equation}
\label{wave_trunc}
\psi^0_k(t)=\sigma_k(t) + \gamma_k\, (t-t^1_{2k+1})_+^{m-1},
\end{equation}
where $\gamma_k\in\RR$ and $\sigma_k$ is a suitable spline from $V_0$.
Applying (\ref{wave_trunc}) to (\ref{two_scale}) we obtain
\begin{equation}
\label{s_trunc2}
s_1(t) = s_0(t) + \sum_k d_k^0 \sigma_k(t) 
         + \sum_k d_k^0\; \gamma_k\; (t-t^1_{2k+1})_+^{m-1}
\end{equation} 
Obviously $\tilde{\sigma}:= s_0 + \sum_k d_k \sigma_k$ is a spline from $V_0$, while 
$\gamma_k\, (t-t^1_{2k+1})_+^{m-1}$ is the contribution from the `new' knot $t^1_{2k+1}$.
By comparison of coefficients in (\ref{s_trunc1}) and
(\ref{s_trunc2}) we conclude that 
$$
\tilde{\sigma}(t) = \sum_{k\in\ZZ} a_{2k} (t-t^1_{2k})_+^{m-1} 
$$
and 
$$
d_k^0 = \frac{a_{2k+1}}{\gamma_k},
$$
i.e., if we know $a_{2k+1}$ and $\gamma_k$ we can easily compute the wavelet
coefficient $d^0_k$.

To determine the coefficients of the truncated power representation
we use the following fact. Obviously the $(m-1)$-th derivative of $s_1$ is given as
$$
s_1^{(m-1)} = (m-1)!\, \sum_{k\in\ZZ} a_k\; \chi_{(t^1_k,\infty)}.
$$
On the other hand, applying the differentiation rule  (\ref{spl_der})
we can compute coefficients $\beta_k$ from the known spline coefficients $c_k^1$
such that
$$
s_1^{(m-1)} =  \sum_{k\in\ZZ}\beta_k\; N^{\cT^1}_{1,k}.
$$
Since $N^{\cT^1}_{1,k}=\chi_{(t^1_k,t^1_{k+1}]}$ we conclude immediately that
\begin{equation}\label{trunc_coeff}
a_k=\frac{\beta_k-\beta_{k-1}}{(m-1)!}.
\end{equation}
Analogously to $a_{2k+1}$ we determine $\gamma_k$ from the
$(m+\widetilde{m}-1)$-th derivative of $\alpha^j_k\,N^{\Xi^0_k}_{m+\widetilde{m},0}$, which is by
(\ref{wave_bit}) the $(m-1)$-th derivative of $\psi^1_k$, i.e., 
\begin{equation}\label{gammak}
\gamma_k = \frac{\rho_{k,\ell_1}-\rho_{k,\ell_1-1}}{(m-1)!},
\end{equation}
where the $\rho_{k,\ell}$ are given by
$$
\alpha^j_k\,\tfrac{d^{m+\widetilde{m}-1}}{dt^{m+\widetilde{m}-1}}N^{\Xi^0_k}_{m+\widetilde{m},0} (t)
=  \sum_{\ell=0}^{m+\widetilde{m}-1}\rho_{k,\ell}\, N^{\Xi^0_k}_{1,\ell}(t).
$$

\paragraph{Coarse scale approximation}
Furthermore, one has to compute the coarse scale approximation 
$$
s_0= \sum_{k\in\ZZ} c^0_k \varphi^0_k = s_1 - \sum_{k\in\ZZ} d^0_k \psi^0_k.
$$

In order to subtract the wavelets from the spline $s_1$ we need the representation
of $\psi_k^0$ in terms of the $\varphi^1_k$. Applying the differentiation
formula (\ref{spl_der}) $\widetilde{m}$ times to $\alpha_k\,N^{\Xi^0_k}_{m+\widetilde{m},0}$ 
(cf.~\ref{wave_bit}) will
yield the representation\footnote{This can be done as part of the computation of $\gamma_k$ above.}
\begin{equation}\label{ukl}
\psi^0_k = \sum_{\ell=0}^{\tilde{m}} u_{k,\ell} N_{m,\ell}^{\Xi^0_k}
\end{equation}
in terms of B-splines over $\Xi^0_k\subset \cT^1$. 
The required representation 
\begin{equation}\label{wave_refin}
\psi^0_k = \sum_{\ell\in\ZZ} b_{k,\ell} \underbrace{N_{m,\ell}^{\cT^1}}_{\varphi^1_\ell}
\end{equation}
is then obtained
by the Oslo algorithm (\ref{oslo}). 
Due to the compact support of $\psi^1_k$ there are only $2\tilde{m}+m-2$ non-vanishing coefficients $b_{k,\ell}$
for each $k$. 
Now, we are able to subtract the wavelet representation
from $s_1$ and obtain
$$
s_0 = \sum_{k\in\ZZ} \tilde{c}_k \,\varphi^1_k
$$
with 
\begin{equation} \label{deco_coarse}
\tilde{c}_k = c^0_k - \sum_{\ell=0}^{n-1} b_{\ell,k}\, d_\ell^0.
\end{equation}

\paragraph{Knot removal}
The above computations yield a representation of $s_0 \in V_0$ in terms of the B-splines
$\varphi_1=N^{\cT^1}_{m,k}$. Since we are looking for a representation in terms of the
$\varphi_0=N^{\cT^0}_{m,k}$, we need a method to
`remove' the knots $t^0_{2k+1}$.
This can be done using an inversion of Boehm's algorithm (\ref{boehm})
(cf.~\cite{EcHa95}). We will remove the knots step by step. This means that
in the $k$-th step we have to replace the representation
$s_1 = \sum_\ell p^k_\ell N_{m,\ell}^{\Theta_k}$ by
$s_1 = \sum_\ell p^{k+1}_\ell N_{m,\ell}^{\Theta_{k+1}}$,
where $\Theta_k=\big(\theta_i^k\big)_{i\in\ZZ}$ with
\begin{equation}\label{theta}
\theta^k_{k+i} =\left\{\begin{array}{ll}
t^0_{k+i}, & \mbox{if }i\le 0\\
t^1_{2k+i},& \mbox{if }i> 0,
\end{array}\right. 
\end{equation}
i.e., $\cT^1=\Theta_0$ and $\cT^0=\Theta_n$.
Starting from $p^0_\ell=\tilde{c}_\ell$ we obtain the required coefficients as
$c^1_\ell = p^n_\ell$.

Following Boehm's algorithm in (\ref{boehm}) the coefficients are related by
$$
p^k_i=\left\{\begin{array}{ll}
p^{k+1}_i,&\mbox{if }i\le k-m\\[1ex]
\displaystyle\frac{(\theta^{k+1}_{i+m}-t_{2k+1})\, p^{k+1}_{i-1}
               +(t_{2k+1}-\theta^{k+1}_i)\,p^{k+1}_i}{\theta^{k+1}_{i+m}-\theta^{k+1}_i},&
               \mbox{if }k-m<i< k,\\[1ex]
p^{k+1}_{i-1},&\mbox{if }i\ge k.
\end{array}\right.
$$
Thus, given $p_i^k$ we obtain $p^{k+1}_i$ by $p^{k+1}_i=p^k_i$, $i\le k-m$,
$p^{k+1}_i=p^k_{i+1}$, $i\ge k-1$, and either
\begin{equation}\label{ib_forward}
p^{k+1}_i = \frac{ (\theta^{k+1}_{i+m}-\theta^{k+1}_i)\, p^k_i
               -(\theta^{k+1}_{i+m}-t_{2k+1})\, p^{k+1}_{i-1}}{t_{2k+1}-\theta^{k+1}_i}
\end{equation}
or
\begin{equation}\label{ib_backward}
p^{k+1}_i = \frac{ (\theta^{k+1}_{i+1+m}-\theta^{k+1}_{i+1})\, p^k_{i+1}
               -(t_{2k+1}-\theta^{k+1}_{i+1})\, p^{k+1}_{i+1}}{\theta^{k+1}_{i+1+m}-t_{2k+1}}
\end{equation}
for $k-m<i<k-1$.
Both recursions yield the same result, if computed in exact arithmetic, since
$s_1\in S_{m,\cT^1}$ and $t_{2k+1}\notin \cT^1$. 

However, rounding in
floating point computations may lead to different results.
Numerical test show that (\ref{ib_forward}) is numerically unstable, while 
(\ref{ib_backward}) yields the expected result. This behavior can be explained 
by the following observation.
Obviously $c_{k-m}=p^k_{k-m}$, i.e., $c_{k-m}$ is determined after
the removal of $k-1$ knots. In (\ref{ib_forward}) the computation of
$p^{k+1}_{k-m+1}$ involves $p^{k}_{k-m}=p^{k+1}_{k-m}$, which in turn depends on
computations of the previous steps down to $p^1_{1-m}$. Thus 
rounding errors of all computations can accumulate and the error may increase with $k$.
In practice it can indeed be observed that the error exceeds the magnitude of the coefficients
for moderate sized $k$.

On the other hand, for (\ref{ib_backward}) the coefficient $p^k_{k-m}$ depends only
on $p^{k-1}_i$, $i=k-m+1,\ldots,k-1$. Successive application of this
argument yields that $p^k_{k-m}$ depends only on $p^{k-j}_i$,
$i=k-m+j,\ldots,k-2+j$. Since $p^{k-m}_i =\tilde{c}_{i+k-m-1}$,
$i=k,\ldots,k-2+m$, it follows that rounding errors made until step $k-m-1$
have no influence on result for $p^k_{k-m}$. Thus, the accuracy of the floating point computations
does not depend on $k$ and the
computation by (\ref{ib_backward}) is stable, which is confirmed by numerical tests. 

The difference in stability is due to the fact, that we are removing the
knots from
`left' to `right' (cf. (\ref{theta})). The recursion (\ref{ib_forward}), which is 
going also from left
to right starts from $p^{k+1}_{k-m}=p^k_{k-m}$, which contains the rounding error
from the previous knot removal. The recursion (\ref{ib_backward}) goes
from right to left starting from $p^{k+1}_{k-1}=p^k_k$, which was not touched by
previous computations.

Summing up we have the following algorithm for the wavelet decomposition.

\algo{\label{deco}
\textbf{Wavelet decomposition}\\[1ex]
\begin{tabular}{lll}
\textsc{Inputs:} & $m$       & spline order\\
                 & $\widetilde{m}$      & number of vanishing moments for wavelets\\
                 & $\big(t^0_i\big)_i$ & spline knot sequence\\
                 & $\big(c^1_i\big)_i$ & vector of spline coefficients\\[0.5ex]
\textsc{Outputs:}& $\big(c^0_i\big)_i$ & vector of spline coefficients\\
                 & $\big(d^0_i\big)_i$ & vector of wavelet coefficients\\[0.5ex]
\textsc{Code:}
\end{tabular}

\setlength{\leftmargini}{30pt}
\setlength{\leftmarginii}{15pt}
\setlength{\leftmarginiii}{15pt}
\begin{enumerate}
\item
{\it Wavelet coefficients and coarse scale approximation}
\textsf{
\begin{enumerate}
\item 
Compute truncated power coefficients $(a_{2k+1})_k$ by \textrm{(\ref{spl_der})} 
and \textrm{(\ref{trunc_coeff})}
\item
FOR $k=0,\ldots,n-1$\quad {\it (loop over wavelets)}
\begin{enumerate}
\item
Compute coefficients $b_{\ell,k}$ for \textrm{(\ref{wave_refin})} by \textrm{(\ref{spl_der})} and Oslo algorithm
\textrm{(\ref{oslo})}
\item 
Compute truncated power coefficient $\gamma_k$ by \textrm{(\ref{spl_der})} and \textrm{(\ref{gammak})}
\item
$d_k^0= \frac{a_{2k+1}}{\gamma_k}$
\item
FOR $\ell = 2(k+1-\ell_1),\ldots,2(k+\ell_2)-m$\\[0.5ex]
\hspace*{2em} $c_\ell^1 := c_\ell^1 - b_{k,\ell}\, d_k^0 $ 
\qquad{\it (subtract $d^0_\ell \psi^0_\ell$ from $s_1$)}
\end{enumerate}
\end{enumerate}}
\item
{\it Knot removal}
\textsf{
\begin{enumerate}
\item
$p_i^0=c_i^1$, $i=0,\ldots,n$
\item
FOR $k=0,\ldots,n-1$ \quad{\it (loop over wavelet knots)}\\[0.5ex]
\hspace*{2em} Compute $p_i^{k+1}$ by \textrm{(\ref{ib_backward})}
\item 
$c_i^0 := p_i^n$
\end{enumerate}}
\end{enumerate}
}

Obviously the number of floating point operation is bounded by $C\, n$, with some constant depending
on $m$ and $\widetilde{m}$, i.e., we have indeed presented a fast algorithm. However, $C$ is considerably larger
compared to the uniform case. 

For equidistant knots, we have due to translation invariance that $c_{k,\ell}^j=p_{2\ell-k}$ and 
$d_{k,\ell}^j=q_{2\ell-k}$. Since the $m+2\widetilde{m}$ non-vanishing coefficients $p_\ell$ and
the $m$ non-vanishing coefficients $q_\ell$ can be determined in advance, the well known classical
algorithm based on the decomposition relation (\ref{dekrel}) will need $4(m+\widetilde{m})n$ floating point operations.

The performance of Algo.~\ref{deco} is essentially reduced by the computations in steps 1(b)i--iii, which depend
only on the used wavelets, but not on the input spline. If we want to decompose several splines over the same
grid, these computations are repeated again and again. To avoid this effect, we collect all splines into a vector valued
function, which means that the corresponding coefficients $c^j_k$, $d^j_k$, $a_k$, and $p^k_i$ in the algorithm
become vectors from $\RR^N$. 
This approach avoids extra memory to store quantities as $b_{l,k}$ and $\gamma_k$, and
is motivated by our applications in circuit simulation as described in \S\ref{num_test}.  
Counting the floating point operations for each step yields then
\smallskip

\begin{tabular}{lrr}
1(a): && $(4N+4)(m-1)n+n\,N$\\
1(b)i: (\ref{oslo})&& $\frac{7}{2}(2\widetilde{m}+m-1)(m-1)m\,n$\\
1(b)ii and 2(b)i (\ref{spl_der}):\qquad & $n+n\sum_{k=0}^{\widetilde{m}+m-1}4 k=$
                                         & $n+ 2(\widetilde{m}+m)(\widetilde{m}+m-1)\, n$\\
1(b)iii:  && $n$\\
1(b)iv: &&$2(2\widetilde{m}+m)n\,N$\\
2.&& $(4N+3)(m-2)n$\\\hline
Total:&&$N(4\widetilde{m}+10m-11)n$\\
    \multicolumn{3}{r}
{$+\big\{m^2(3.5m + 7\widetilde{m}-5)
           +\widetilde{m}(2\widetilde{m}-3m-2)+8.5m-8\big\}n$}\\\hline\hline
\end{tabular} 
\medskip

Here, we have used that the computation of the $u_{k,\ell}$ in (\ref{ukl}) is an intermediate result in the computation 
of $\gamma_k$, while the division by $(m-1)!$ in (\ref{trunc_coeff}) and (\ref{gammak}) can be omitted yielding the 
same result. Obviously, for sufficiently large $N$ the algorithm will need two to three times as much operations as
the classical wavelet transform on a grid of the same size. However, a nonuniform grid may be chosen much smaller in some 
cases, due to its flexibility. Thus, the use of nonuniform spline wavelets may be beneficial for suitable applications, as we 
will demonstrate in an example in \S\ref{num_test}.  

The decomposition algorithm can be applied again to $s_1$, which is decomposed in a
coarser signal $s_2$ and details $w_2$. The successive application of this approach,
also known as the pyramid scheme, yields the multiscale decomposition
\begin{equation}\label{multi_scale}
s_J=s_0+ \sum_{j=0}^{J-1} \sum_{k\in\ZZ} d^j_k \psi^j_k.
\end{equation}

\subsection{Fast wavelet reconstruction\label{wv_reko}}

In order to transform a wavelet representation into the corresponding
B-spline representation, we have to solve the following reverse problem.
Given coefficients $c^0_k$ and $d^0_k$, we have to compute coefficients
$c^1_k$ which satisfy (\ref{two_scale}).
In fact this problem can be solved immediately by the Oslo algorithm yielding
$s_0 = \sum_{k\in\ZZ} \tilde{c}^1_k \varphi^1_k$
and $\psi^1_k = \sum_{\ell\in\ZZ} b_{k,\ell} \varphi^1_\ell$.
Then we obtain the required coefficients immediately as
$c^1_k = \tilde{c}^1_k + \sum_{\ell\in\ZZ} b_{\ell,k} d_\ell^0$.

\algo{\label{reco}
\textbf{Wavelet reconstruction}\\[1ex]
\begin{tabular}{lll}
\textsc{Inputs:}& $m$        & spline order\\
                & $\widetilde{m}$      & number of vanishing moments for wavelets\\
                & $\big(t^0_i\big)_i$ & spline knot sequence\\
                & $\big(c^0_i\big)_i$   & vector of spline coefficients\\
                & $\big(d^0_i\big)_i$   & vector of wavelet coefficients\\[0.5ex]
\textsc{Output:}& $\big(c^1_i\big)_i$   & vector of spline coefficients\\[0.5ex]
\textsc{Code:}
\end{tabular}
\textsf{
\setlength{\leftmargini}{30pt}
\setlength{\leftmarginii}{15pt}
\setlength{\leftmarginiii}{15pt}
\begin{enumerate}
\item 
Compute $\big(\tilde{c}^1_i\big)_i$ by the Oslo algorithm  \textrm{(\ref{oslo})}
\item
$\big(c^1_i\big)_i=\big(\tilde{c}^1_i\big)_i$
\item
FOR $k=0,\ldots,n-1$ \quad{\it (loop over wavelets)}
\begin{enumerate}
\item
Compute coefficients $b_{k,\ell}$ in \textrm{(\ref{wave_refin})} by successive
differentiation \textrm{(\ref{spl_der})} and Oslo algorithm 
\textrm{(\ref{oslo})}
\item 
FOR $\ell = 2(k+1-\ell_1),\ldots,2(k+\ell_2)-m$\\[0.5ex]
\hspace*{2em} $c_\ell^1 := c_\ell^1 + b_{k,\ell}\,d_k^0$ \qquad
{\it (Add $d^0_\ell \psi^0_\ell$ to $s_1$)}
\end{enumerate}
\end{enumerate}
}}

Analogously to Algo.~\ref{deco} we determine the computational cost of Algo.~\ref{reco} for vector valued input

\begin{tabular}{lrr}
1: &&$(4N+3)m(m-1)n$\\
3(a) (\ref{spl_der}): & \hspace{8em}$n\sum_{k=0}^{\widetilde{m}}4 k=$
                                         & $2\widetilde{m}(\widetilde{m}+1)\, n$\\
3(a) (\ref{oslo}):&& $\frac{7}{2}(2\widetilde{m}+m-1)(m-1)m\,n$\\
3(b): &&$2(2\widetilde{m}+m)n\,N$\\\hline
Total:&&$N\big(4\widetilde{m}+m(4m-2)\big)n$\\
  &  \multicolumn{2}{r}
{$+\big\{m^2(3.5m + 7\widetilde{m}-4)
           +\widetilde{m}(2\widetilde{m}-7m+6)+2.5m\big\}n$}\\\hline\hline
\end{tabular} 
\medskip

Again the algorithm is more costly than the classical algorithm for a uniform grid with $4N(m+\widetilde{m})n$ operations.
However, if a much smaller, adapted grid can be used a nonuniform grid is of advantage.
Another question is how a suitable nonuniform grid can be determined. Questions like this have lead to the
investigation of non-linear approximation methods. For a survey we refer the reader to \cite{DeV}.
In Sect~\ref{wave_adapt}
we present methods to determine an optimal adaptive grid using modifications of our wavelet algorithms.

\section{Modifications of the wavelet definition\label{modify}}

\subsection{Spline wavelets on the interval}

In many cases, one needs spline representation on a compact interval $[a,b]$. 
It is reasonable to choose the end points as spline knots, i.e., $t_0 = a$ and $t_n = b$.
and define the spline space
$$
S_{m,\cT}([a,b]):=\Big\{s\big|_{[a,b]}:~s\in S_{m,\cT}\Big\},
$$ 
i.e., the restriction of $S_{m,\cT}$ onto $[a,b]$. Obviously, 
$$
S_{m,\cT}([a,b]):=\spn\Big\{N_{m,k}\big|_{[a,b]}:~k=1-m,\ldots,n-1\Big\},
$$ 
and only the knots $t_{1-m},\ldots,t_{n+m-1}$ play a role for the spline space and
the B-splines. 

The spline space $S_{m,\cT}([a,b])$ does not depend on the particular choice of the outer knots
$t_{-i}$ and $t_{n+i}$, $i>0$, but they influence the B-Splines at the boundary.
For stability reasons and simplicity one chooses multiple knots at the boundary, i.e.,
$t_{-i}=t_0= a$ and $t_{n+i} = t_n = b$, $i>0$. 

Wavelets can be defined in principle as before with nested knot sequences $\cT^j=(t^j_k)_{k\in\ZZ}$, $j\in\No$,
$$
t^j_k=\left\{\begin{array}{ll}
a,&\mbox{if } k\le 0,\\
b,&\mbox{if } k\ge n_j,\\
t^{j+1}_{2k}&\mbox{otherwise,}
\end{array}\right.
$$
where $n_j=\big\lceil n_{j+1}/2\big\rceil$
However, to ensure $\tilde{m}$ vanishing moments
we need all derivatives of  $N^{\Xi^j_k}_{m+\widetilde{m},0}$ to vanish at the boundaries.
Therefore, we modify the definition of $\Xi^j_k$ in (\ref{wave_knots}) to
$$
\Xi^j_k:=\left\{\!\!\!\!\begin{array}{lr}
\big(t^j_{2-m},t^j_{3-m},\ldots,
   t^j_k,t^{j+1}_{2k+1},t^j_{k+1},\ldots,
   t^j_{\tilde{m}+1}\big),&\hspace{-1em}\mbox{if }0\le k<1-m+\ell_1\\[.5ex]
\big(t^j_{n_j-\tilde{m}-1},..,
   t^j_k,t^{j+1}_{2k+1},t^j_{k+1},..,
   t^j_{n_j+m-2}\big),&\hspace{-6.8em}\mbox{if }n_j+m-2-\ell_2<k<n_{j+1}-n_j,\\[.5ex]
\big(t^j_{k+1-\ell_1},t^j_{k+2-\ell_1},\ldots
   t^j_k,t^{j+1}_{2k+1},t^j_{k+1},\ldots,
   t^j_{k-1+\ell_2},t^j_{k+\ell_2}\big)&\mbox{otherwise.}
\end{array}\!\!\!\!\right.
$$
With these settings the above algorithms can be applied as before, using the modification 
for multiple knots described in Sect~\ref{splines}.

\subsection{Periodic spline wavelets}

Many problems deal with periodic functions. Our approach can easily be modified to define 
periodic spline wavelets. A $P$-periodic spline space can be defined if the knots satisfy the periodicity
condition
$$
t_{k+n}=t_k+P,\qquad k\in\ZZ
$$
for some $n\in\NN$. Then the periodic spline space can be defined as
\begin{eqnarray*}
\nonumber
S^P_{m,\cT}&:=&\Big\{s\in S_{m,\cT}:~s(t)=s(t+P),~t\in\RR\Big\}\\
&=&\Big\{\sum_{k\in\ZZ}c_k N_{m,k}^\cT:~c_k=c_{k+n}\Big\}
=\spn\left\{ N_{m,k}^{\cT,P}:~k=1,\ldots,n\right\}
\end{eqnarray*}
with the periodic B-splines 
$$
N_{m,k}^P:=N_{m,k}^{\cT,P}:=
\sum_{\ell\in\ZZ} N_{m,k}^\cT(\cdot+\ell\,P)=\sum_{\ell\in\ZZ} N_{m,k+\ell\,n}^\cT.
$$

To define wavelets we introduce again nested knot sets $\cT^j=(t^j_k)_{k\in\ZZ}$, $j\in\No$,
satisfying $t^j_k+P =t^j_{k+n_j}$, $k\in\ZZ$, $n_j=\big\lceil n_{j+1}/2\big\rceil$, 
and $t^j_k=t^{j+1}_{2k}$, $k=0,\ldots,n_j-1$.
The corresponding periodic spline wavelets are then defined as
$$
\psi^P_{j,k}:=\sum_{\ell\in\ZZ} \psi^j_k(\cdot +\ell\,P)=\sum_{\ell\in\ZZ} \psi^j_{k+\ell(n_{j+1}-n_j)}.
$$

Now the algorithms introduced above can be immediately adapted for periodic splines, using the fact
that any $P$-periodic spline is uniquely determined by the vectors $(t_k)_{k=0}^{n-1}$ and  $(c_k)_{k=0}^{n-1}$.
All of the above algorithms can be applied using $t_k=t_{\kappa_1(k)}+\kappa_2(k)\,P$ for spline knots and
$c_k=c_{\kappa_1(k)}$ for coefficients, where the mapping $\kappa:\ZZ\to\{0,\ldots,n-1\}\times\ZZ$ is uniquely
determined by $k=\kappa_1(k)+n\,\kappa_2(k)$. 

\subsection{Further generalizations}

Following the ideas of \cite{Bit05a,Bit05b} and \S\ref{wavelets}--\ref{algorithms} one can consider further
generalization. However, this leads to a more complex notation, while the principal ideas stay the same.
Since the goal of this publication is to present an efficient and plain implementation with easily maintainable code
we did not consider all possible generalizations, but only those which may give an advantage in our target applications. 
However, since the focus may change with other problems under consideration, we will give here 
some suggestions on possible generalized settings, where the ideas from our approach can still be applied.

\paragraph{Choice of $\Xi^j_k$} Apparently we could replace $\ell_1$ and $\ell_2$ in (\ref{wave_knots}) 
by any $\ell_{1}(k)$ and $\ell_{2}(k)$,
as long  as $\ell_{1}(k)+\ell_{2}(k)=m+\tilde{m}$. Such a modification has been done for the boundary wavelets on
the interval. However, the current choice was made to get some symmetry, which is often of advantage. 

More general settings for the knot sequence $\Xi^j_k$ are possible as long as it contains  $t^{j+1}_k$ and $m+\tilde{m}$
knots from $\cT^j$. The algorithms would follow the same principle, but the support size would increase
and the implementation would become more involved
with no obvious advantage.

\paragraph{Choice of grids $\cT^j$} As long as $\cT^j\subset\cT^{j+1}$, any sequence of knot sets can be used.
For each new knot $t_k^{j+1}\not\in\cT^j$ a wavelet is defined (\ref{wave_bit}), where $\Xi^j_k$  contains  $t^{j+1}_k$ 
and $m+\tilde{m}$ additional knots from $\cT^j$ (which should be close to $t^{j+1}_k$). The algorithms would work analogously.
However, one has to deal with the more complex setting, while the handling of the grids would require extra memory and computation time.

\section{Wavelets for grid adaptation\label{wave_adapt}}

An interesting property of wavelet is that they permit an efficient adaptive approximation
using the best $n$-term approximation. That is, the function is approximated by a linear combination
of $n$ adaptively chosen wavelets. For a stable wavelet basis the wavelets with the largest expansion coefficients
yield the best approximation. The approximation is optimal for functions from Besov spaces (see e.g.~\cite{DeV}).
The same approximation power is obtained with free knot splines, i.e., B-spline representations
with adaptively chosen knots.
This is apparent, since any spline wavelet representation is a particular spline
with adaptive knots (but not vice versa). 

A simple rule of thump is that an adaptive spline 
scheme is of benefit for functions with isolated singularities (e.g. discontinuities of
the function or its (higher order) derivatives, sharp transients). By a local refinement, i.e., additional wavelets or additional spline
knots in the vicinity of such isolated singularities, one can achieve an accurate approximation with a relative small number
of degrees of freedom. Going for the optimum means to determine the free knot spline of \emph{best approximation}, which
minimizes the approximation error for a given number of spline knots. Often it is not possible to find this optimal 
approximation at reasonable cost such that one looks instead for an \emph{almost} best approximation, where the error has
at least the same order of magnitude as the best approximation.

In this section we show how our spline wavelets can be used to generate 
an adaptive grid of spline knots. To achieve this goal we will use the fact that the wavelet $\psi^j_k$ corresponds 
by definition with the knot $t^{j+1}_{2k+1}$.

\subsection{Coarsening\label{coarse}}
In practice one can obtain an efficient $n$-term approximation as follows. We start from a sufficiently accurate
spline approximation, i.e., an expansion in terms of scaling functions $\varphi^J_k$ with sufficiently large $J$,
which can be obtained e.g. by interpolation or quasi-interpolation. Next the fast decomposition algorithm is used
to compute the multiscale wavelet expansion (\ref{multi_scale}). From this expansion the terms with small 
coefficients are
removed (thresholding), which may reduce the size of the expansion essentially, while only a small error is introduced.

Let us consider the simple case that $J=1$.
That is an application of the decomposition algorithm yields
\begin{equation}\label{deco1}
s_1=s_0+\sum_k d^0_k \psi^0_k.
\end{equation}
Thresholding means that one replaces $s_1$ by   
\begin{equation}\label{decothresh}
\tilde{s}_1=s_0 + \sum_{k:~|d^0_k|\ge\varepsilon} d^0_k \psi^0_k.
\end{equation}
Subtracting (\ref{decothresh}) from (\ref{deco1}) one obtains
$$
\tilde{s}_1=s_1 - \sum_{k:~|d^1_k|<\varepsilon} d^1_k \psi^1_k.
$$
Obviously $\tilde{s}_1$ can be computed by a modified version of the decomposition algorithm.
In particular, Step 1(b)iv and 2(b) (with correspondingly modified $\Theta_k$) in Algorithm~\ref{deco} 
are only performed if  $|d^1_k|<\varepsilon$.  
The result is $\tilde{s}_1$ in the B-spline expansion over the coarsened grid
$$
\tilde{\cT}_1:=\cT^1\setminus\{t^1_{2k+1}:~|d^1_k|<\varepsilon\}.
$$
The above method can be applied successively to
$\tilde{s}_1$ in order to obtain further coarsened representations
$\tilde{s}_j$ on $\tilde{\cT}_j$. 

In this process the spline grid will be adapted according to the approximated function. Starting from an approximation
with required accuracy on a uniform or nonuniform 
grid unnecessary knots will be removed. This will mainly happen in areas, where the function is smooth such that it can
be approximated by a polynomial of degree $<m$. In localities where the function is less smooth (e.g.\ isolated singularities
or fast transients) few or no knots will be removed to keep the error under control.

The error can be controlled as follows. The normalization factor $\alpha^0_k$ in (\ref{wave_bit}) is chosen
such that the coefficients $b_{k,\ell}$ in (\ref{wave_refin}) satisfy
\begin{equation}\label{normalization}
\max_\ell \big|b_{k,\ell}\big|=1.
\end{equation}
In practice this is done by computing $b_{k,\ell}$ and $\gamma_k$ for $\alpha^0_k=1$, and divide these
coefficients by  $\max_\ell \big|b_{k,\ell}\big|$ afterwards. Then we have
$$
\big|\psi^0_k(t)\big|<\max_\ell \big|b_{k,\ell}\big|=1.
$$
Due to their compact support (see \S\ref{wavelets}), all but $m+\tilde{m}-1$ wavelets vanish for any $t\in\RR$.
Therefore
$$
\big|\tilde{s}_1(t)-s_1(t)\big|=\bigg| \sum_{k:~|d^1_k|<\varepsilon} d^1_k \psi^1_k(t)\bigg|
\le \varepsilon  \sum_k \big|\psi^1_k(t)\big|\le (m+\tilde{m}-1)\varepsilon
$$
and
$$
\big|\tilde{s}_j(t)-s_1(t)\big|\le (m+\tilde{m}-1)\,j\,\varepsilon.
$$

\subsection{Refinement}
It is not always possible to obtain an approximation on a very fine grid, and if possible
it may be to expensive. In particular, for the approximative solution of operator 
equations the reason of an adaptive grid is to reduce the computational cost for the
(often computational expensive) solver,
by reducing the number of degrees of freedom. That is, the number of degrees of freedom
has to be kept as small as possible. 

Let us assume we are able to compute a spline approximation of a function $f$ for any grid in finite time,
which increases at least linearly with the grid size.
Namely, we denote by  $A_{\cT}(f)$ an approximation method  yielding a
spline approximation $f_{\cT}  = A_{\cT}(f)\in S_{m,\cT}$ of $f$ on the grid ${\cT}= (t_k)_{k=0}^n$.
Typical examples for $A_{\cT}(f)$ are collocation or Galerkin solvers for a differential equation.

Our approach is now to start with a coarse initial grid $\cT^0$ and determine a first approximation
$f_0:=A_{\cT^0}(f)$. Now we want to use information contained in  $f_0$ to generate 
a new improved grid. Performing a wavelet decomposition (Algo.~\ref{deco}) we obtain 
$$
f_0=\tilde{f} + \sum_k d_k\, \psi^0_k,
$$
Obviously, a large $\big| d_k\big|$ implies a large approximation error of the coarse
scale approximation $\tilde{f}=\sum_k c_k\, \varphi^0_k$ 
in the support of $\psi^0_k$, i.e., in a neighborhood 
of $t_{2k+1}$. Moreover, a large wavelet coefficient may also indicate a large approximation error of  
$f_0$ near $t_{2k+1}$. This assumption is justified for a large family of functions, e.g.\  solutions of typical differential or integral equations. Our goal is to improve the approximation locally in regions with a large local approximation error
by refining the grid locally in this regions. Namely, we insert
$$
n_k:=\left\lfloor |d_k|\,\frac{\alpha}{\max_\ell|d_\ell|} \right\rfloor
$$
equidistant knots into the intervals $(t_{2k},t_{2k+1}) $ and $(t_{2k+1},t_{2k+2}) $, $k=0,\ldots,n-1$, 
which yields the refined
grid $\cT^1$. Obviously knots are only 
inserted if $\big|d_k\big|$ is not below the threshold $\varepsilon := \frac{\max_\ell|d_\ell|}{\alpha}$ and the number 
of inserted knots
grows with $\big|d_k\big|$. The factor $\alpha>1$ controls the refinement rate. Note that the coarse scale
approximation $\tilde{f}=\sum_k c_k\, \varphi^0_k$ is not needed here, i.e., we can stop the decomposition after 
the computation of the wavelet coefficient in the decomposition algorithm (Algo.~\ref{deco}, Step1(b)iii).

Now we determine the improved approximation $f_1:=A_{\cT^1}(f)$. We cannot expect this approximation to be 
sufficiently accurate, but we have gained additional information for another wavelet refinement. Thus, we repeat 
refinements until we have gained the required accuracy, which results in the following algorithm.

\algo{\label{refine}
\textbf{Wavelet based grid refinement}\\
\begin{tabular}{lll}
\textsc{Inputs:} & $m$       & spline order\\
                 & $\widetilde{m}$      & number of vanishing moments for wavelets\\
                 & ${\cT}^0$ & initial spline knot sequence\\
                 & $\varepsilon$ & error tolerance\\
\textsc{Outputs:}  & ${\cT}^i$ & final spline knot sequence\\
                     & $f_i$ &spline approximation on ${\cT}^i$\\[0.5ex]
\textsc{Code:}
\end{tabular}

\textsf{
\setlength{\leftmargini}{30pt}
\setlength{\leftmarginii}{15pt}
\setlength{\leftmarginiii}{15pt}
\begin{enumerate}
\item
$f_0 := A_{{\cT}^0}(f)$ 
\item
FOR $i =1,2,3\ldots$
\begin{enumerate}
\item
Compute the wavelet coefficients $d_k$, $k=0,\ldots,n-1$, of $f_{i-1}$
\item
Determine ${\cT}^i$ by inserting $n_k:=\left\lfloor |d_k|\,\frac{\alpha}{\max_\ell|d_\ell|} \right\rfloor$
 new knots\\ into $(t^{i-1}_{2k},t^{i-1}_{2k+1}) $ and $(t^{i-1}_{2k+1},t^{i-1}_{2k+2}) $ for $k=0,\ldots,n-1$,
\item
$f_i = A_{{\cT}^i}(f)$
\item
IF $\|f_i - f_{i-1}\|<\varepsilon$ stop with solution $f_i$
\end{enumerate}
\end{enumerate}}
}

The above algorithm is of particular interest, if the approximation method is based on an iterative method, e.g.\ 
Newton's method for a nonlinear problem. Then an initial guess and a stopping criterion is required for the iteration, i.e.,
the approximation method $A_{\cT}(f)$ is replaced by $A_{\cT}(f,f^*,\varepsilon)$, 
where $f^*$ is the initial guess and $\varepsilon$ a tolerance
to stop the iteration. In the above algorithm step 2(c) is replaced by  
$f_i = A_{{\cT}^i}(f,f_{i-1},\varepsilon_i)$. That is we use the already known approximation $f_{i-1}$, where
the spline representation for the new grid ${\cT}^i$ is determined bye the Oslo algorithm. Thus, for larger $i$ we have 
an excellent initial guess, such that only a few Newton steps (often only one) are sufficient. Furthermore, for small
$i$, where many (possibly damped) Newton steps may be necessary  due to a bad initial guess, the grid is much smaller
such that the cost for one step can be neglected. Thus, with a suitable refinement rate $\alpha$, the overall
computational cost may be comparable to the cost of few Newton steps on the final grid. 

The stopping criterion in step 2(d) is based on the assumption that the approximation is essentially improved
by the refinement. This assumption is justified if $\alpha$ is chosen sufficiently large to ensure a non negligible growth
of the grid. Then it can be expected that the approximation error is in the range of $\varepsilon$ when the algorithm stops.
Due to the refinement we expect the grid to be adapted to the solution, i.e., coarser at smooth parts of $f$ and finer
at non smooth parts. However, a final coarsening can sometimes improve the efficiency of the approximation by essentially
reducing the grid size while introducing only an insignificant additional error.

\section{The application of wavelet based grid adaptation in circuit simulation\label{num_test}}

The described methods have been implemented as a C++ library, which can be used in suitable applications.
Here, we present results from an adaptive circuit simulation method. Details about the method 
can be found in \cite{BiDau10a,BiDau10b,BiBra12b,BiBra14b} so that we will give only a short introduction
here. By a Modified Nodal Analysis \cite{HRB75,GFM05} (which is based on Kirchhoff's laws and device modells) 
one obtains the cirquit equations
\begin{equation}\label{mna}
\tfrac{d}{dt}q\big(x(t)\big) + \underbrace{i\big(x(t)\big) + s(t)}_{g(x(t),t)}= 0,
\end{equation}
where $x:\RR\to\RR^N$ denotes the unknown voltages and currents,
$q:\RR^N\to\RR^N$ contains charges and fluxes from capacities and intuctivities, $i:\RR^N\to\RR^N$
describes static contributions from resistors, diodes, transistors, etc., while $s:\RR\to\RR^N$ contains
independent sources. Typically the Jacobian of $q(x)$ is not invertible, i.e., (\ref{mna}) is a system of
differential algebraic equations (DAE).

\subsection{Periodic  steady state computation\label{pss}}

To determine the periodic steady state of the circuit, we have to solve a periodic boundary problem, i.e.,
the solution of (\ref{mna}) has to satisfy $x(t)=x(t+P)$.\footnote{An approach for initial value problems (Transient Analysis) 
using splines on the interval is described in \cite{BiDau10a,BiDau10b}.} 

We discretize (\ref{mna}) by a spline-based Petrov-Galerkin method. That is, the unknown solution is expanded 
as a $P$-periodic spline $x=\sum_{k=1}^n c_k\,N^P_{m,k}$ for unknown coefficients $c_k\in\RR^N$. 
Integrating (\ref{envelope}) over a suitable chosen set of $n$ subintervals $[\tau_{\ell-1},\tau_\ell]$ yields 
the system of $n\times N$ nonlinear equations
\begin{equation}\label{nonlinear}
F_\ell(c_1,\ldots,c_n) := \omega(\tau_k)\Big(q\big(x(\tau_\ell)\big)-q\big(x(\tau_{\ell-1})\big)\Big)
  +\int_{\tau_{\ell-1}}^{\tau_{\ell}} g_k\big(x(t),t\big)\,dt,\quad \ell=1,\ldots,n.
\end{equation}
This system can be solved numerically by applying Newton's method.

Due to their support properties
the use of the B-spline expansion (with free knots) has some advantages compared with the use of a 
wavelet expansion (cf.~\cite{BiUr04,CDD03,DSX00}).
First one has to evaluate $x(t)$ many times. For a B-spline expansion this is much more efficient 
compared to a wavelet expansion. Furthermore, the Jacobian of 
$\big( F_\ell(c_1,\ldots,c_n)\big)_{\ell=1,\ldots,n}$ is much sparser for the B-spline expansion.
Since the computation of the Jacobian and solving the linear system require most of the computation 
time this sparsity is essential for the efficiency of the algorithm.  

However, it remains open how the spline knots have to be chosen for an efficient computation. 
This is in particular important if the 
circuit exhibits digital like behavior, i.e., the solution switches between a finite number of states
with sharp transients. Then a grid adapted to the signal
waveforms in the circuit will yield a sparser representation, which makes the simulation much more efficient. 
Thus, the grid refinement described in 
Algo.~\ref{refine} is used to determine an adaptive grid, 
where $x_i=A_{{\cT}^i}(x,x_{i-1},\varepsilon_i)$ denotes the Newton's method for (\ref{nonlinear}).

For the described circuit simulation problem the advantages of our new algorithms are fully exploited.
Most of the computation time is consumed by the (sparse) linear solver and by the evaluation
of $q$ and $g$ with their Jacobians, which are usually expensive due to the complex models for semiconductor devices. The remaining computations including the wavelet transforms take only $\approx$1\% of the simulation time. Thus, the little extra computation 
time, which may be needed compared to classical wavelet transform, fully pays of due to smaller linear systems and fewer
function evaluations. 

\paragraph{Diode rectifier}

\begin{figure}[htb]
\centerline{\includegraphics[width=0.5\linewidth]{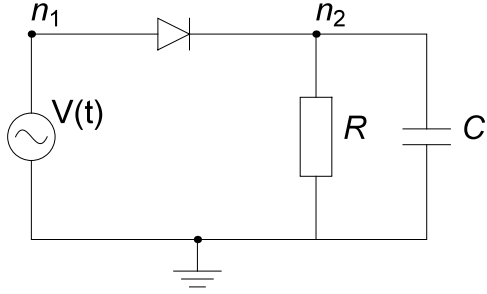}}
\caption{Diode Rectifier \label{rect_scheme}}
\end{figure}

\begin{figure}[htb]
\parbox[t]{0.5\linewidth}{
\includegraphics[width=\linewidth]{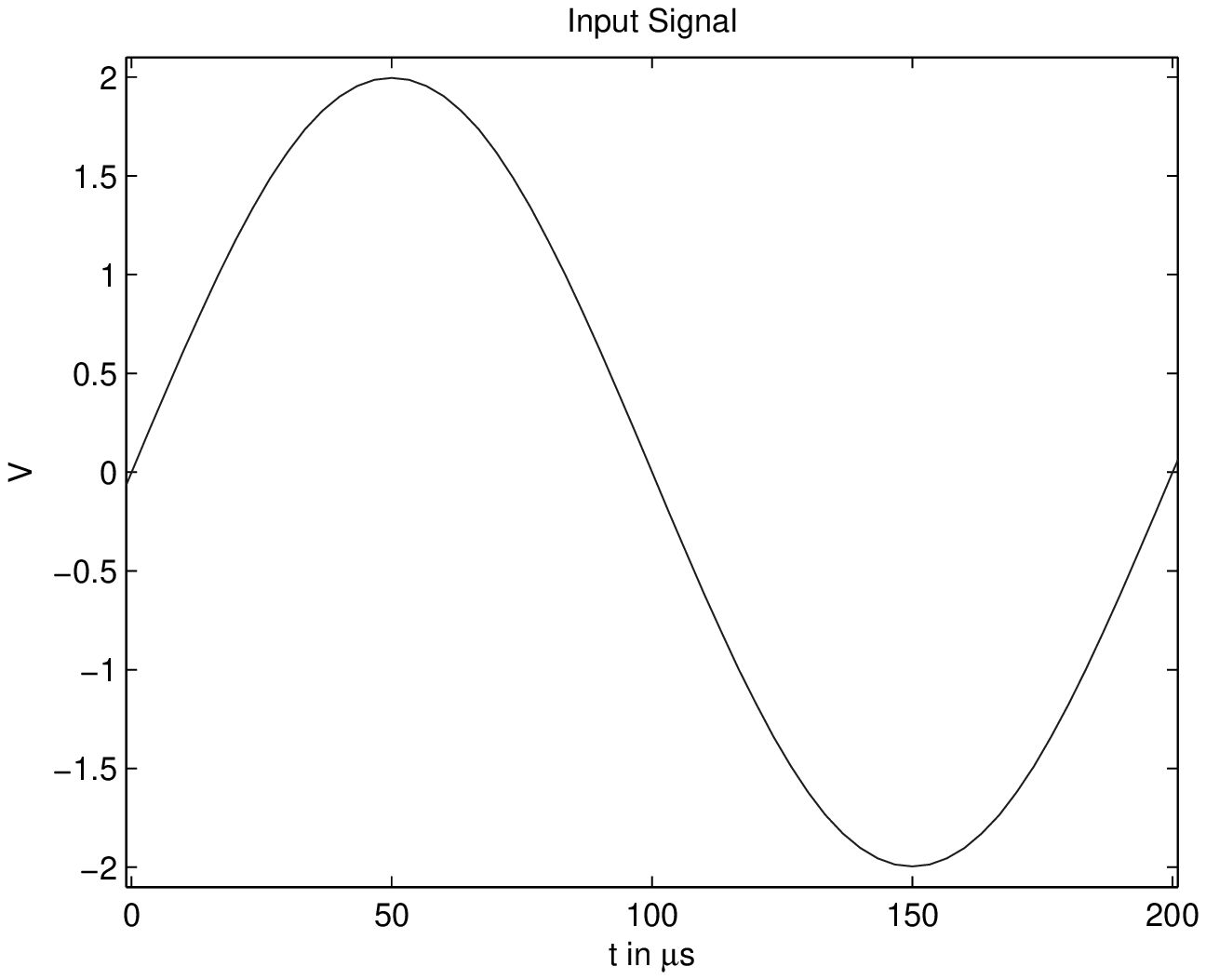}}
\hfill
\parbox[t]{0.5\linewidth}{
\includegraphics[width=\linewidth]{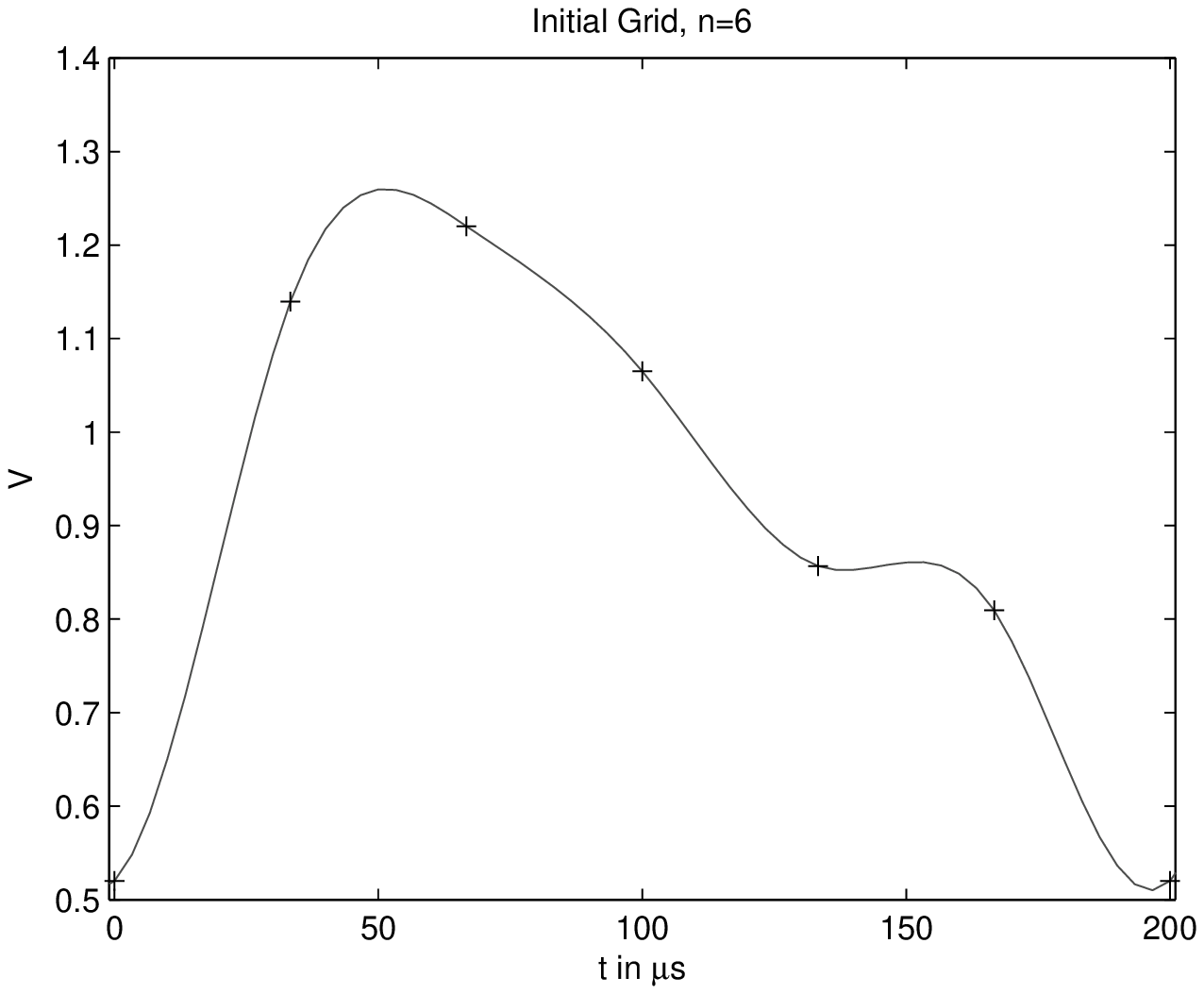}}

\parbox[t]{0.5\linewidth}{
\includegraphics[width=\linewidth]{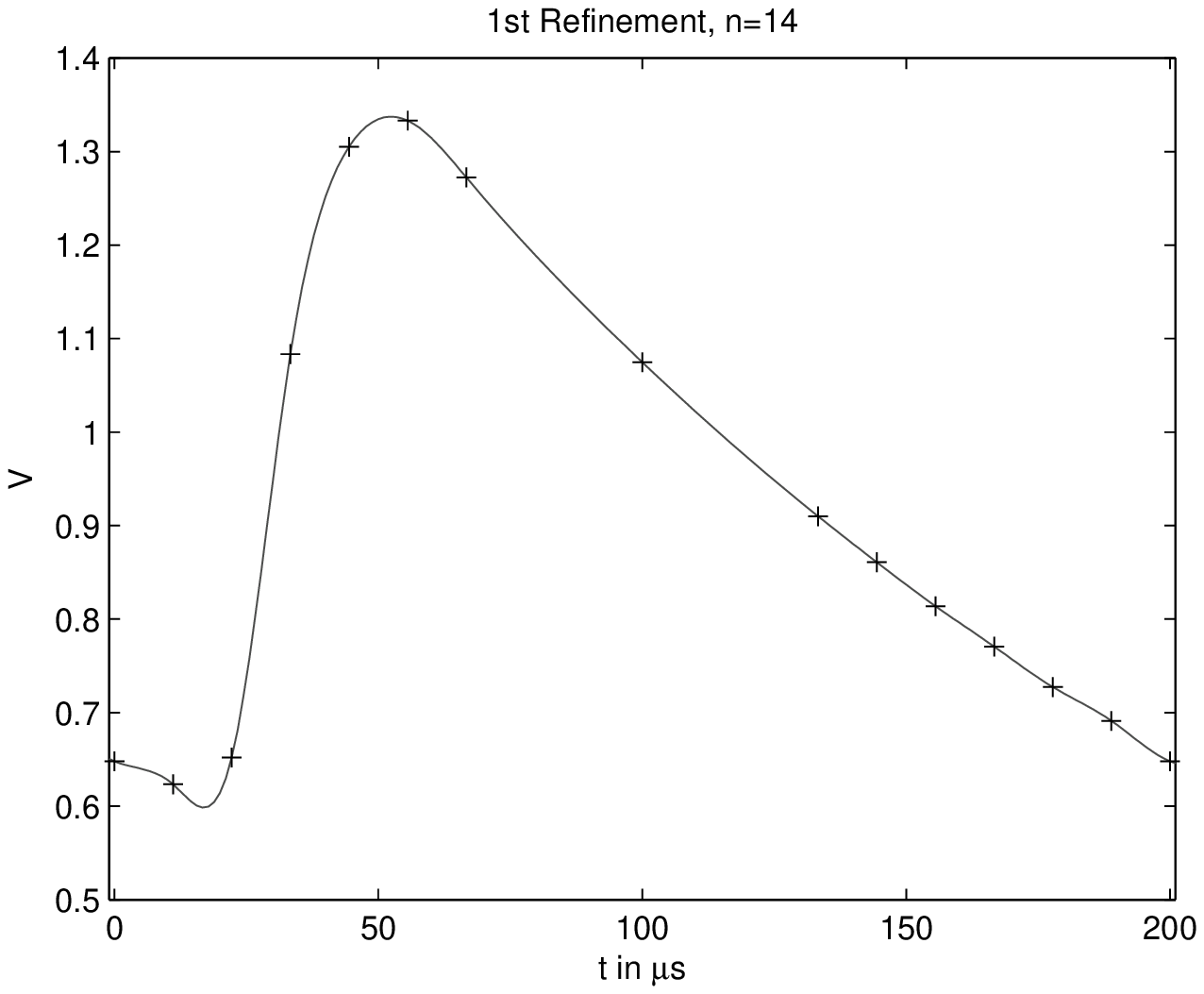}}
\hfill
\parbox[t]{0.5\linewidth}{
\includegraphics[width=\linewidth]{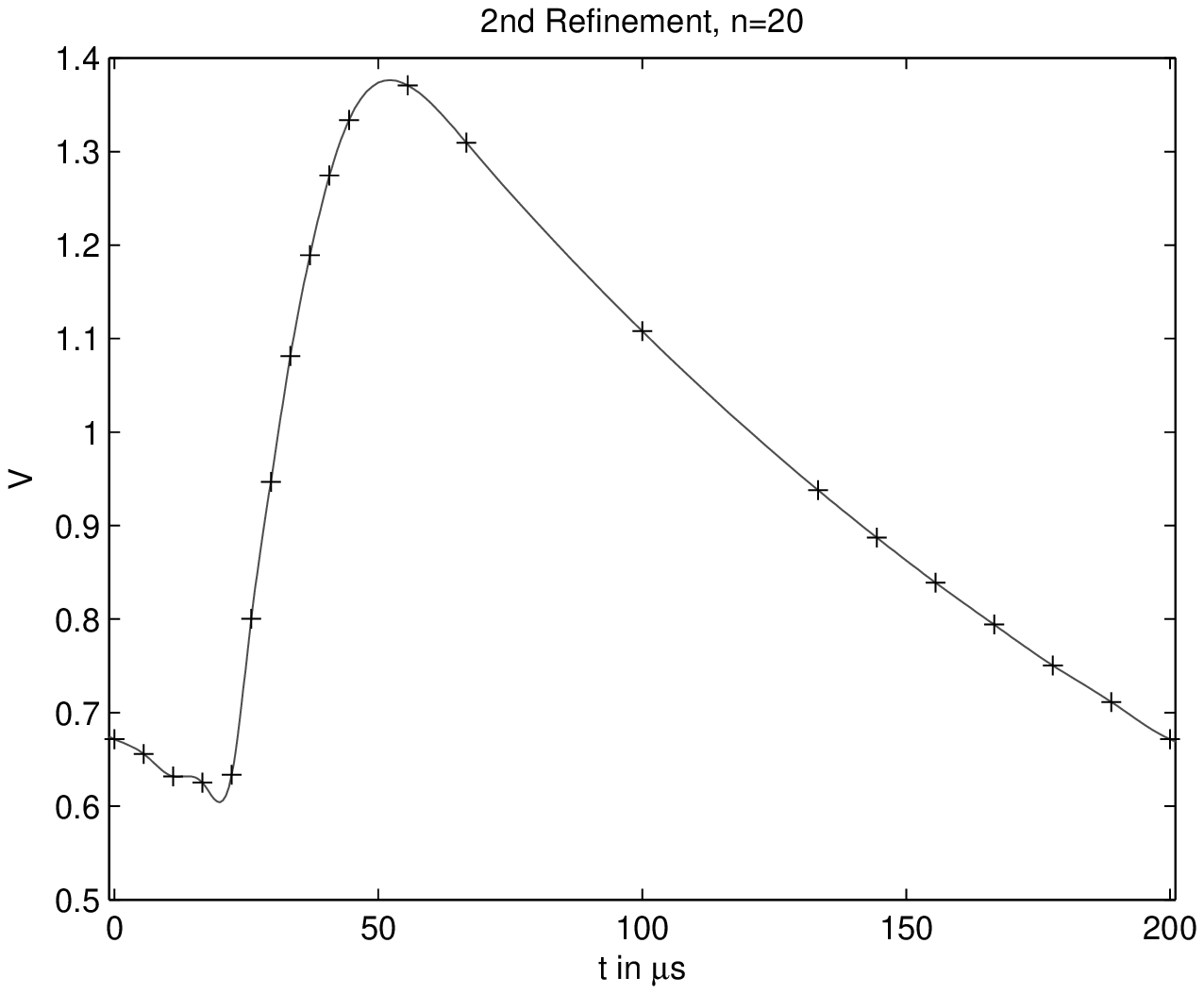}}

\parbox[t]{0.5\linewidth}{
\includegraphics[width=\linewidth]{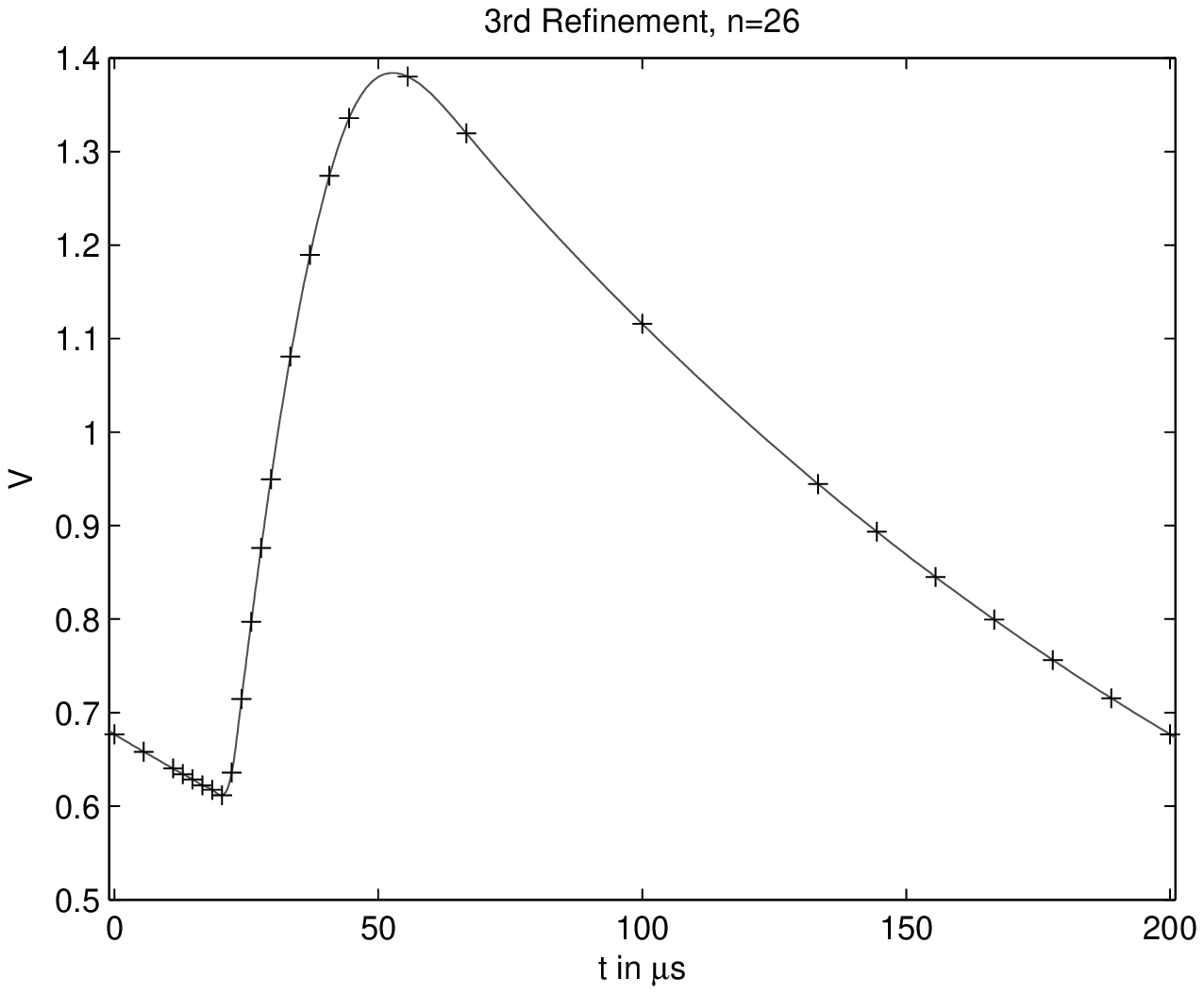}}
\hfill
\parbox[t]{0.5\linewidth}{
\includegraphics[width=\linewidth]{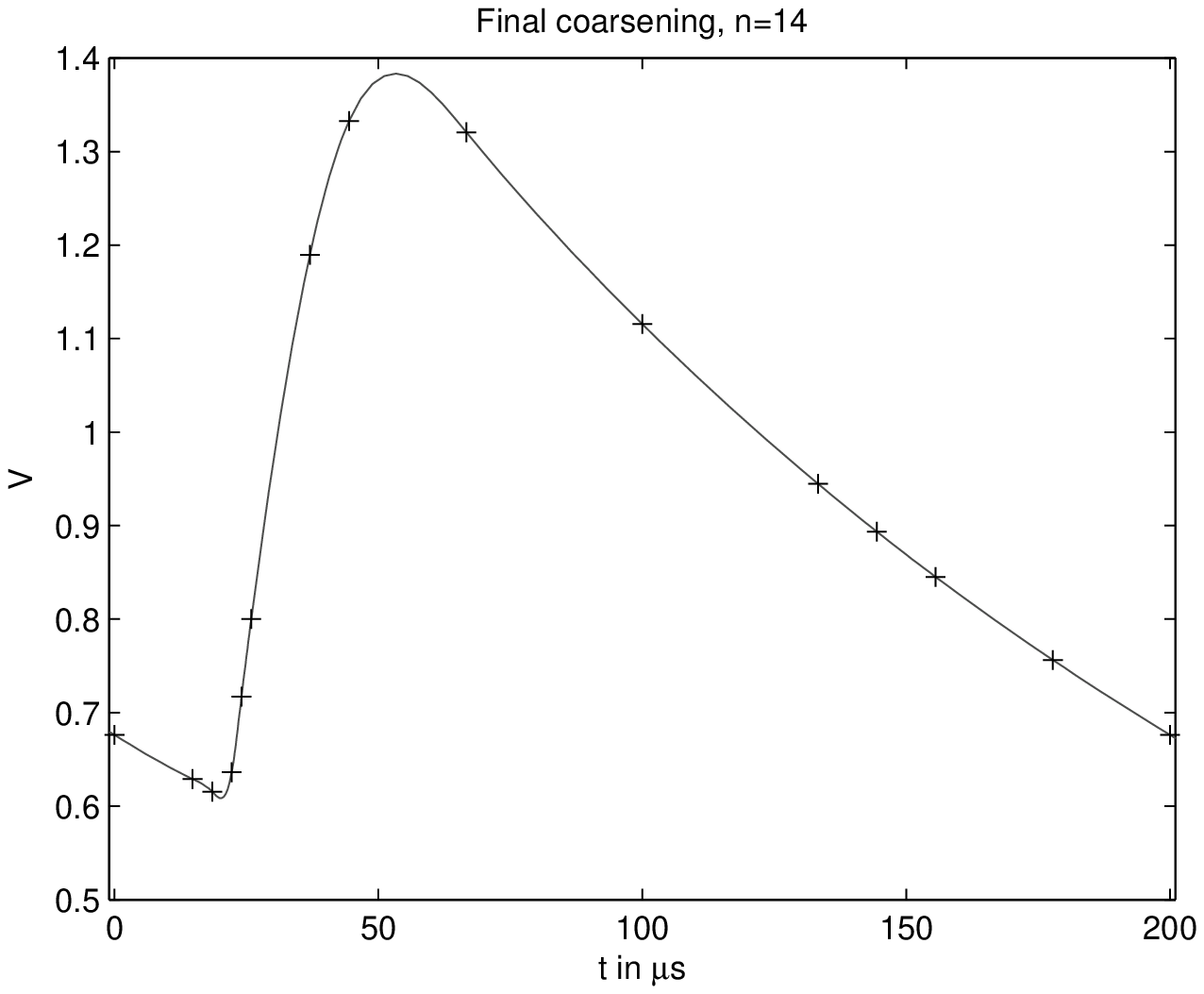}}

\caption{Simulation Results for Diode Rectifier\label{rect}}
\end{figure}

As a first example we show simulation results for a diode rectifier (Fig.~\ref{rect_scheme}). This is a very simple circuit, but 
it permits to demonstrate in illustrative manner how the refinement algorithm works. 
The diode becomes conductive if the branch voltage exceeds
a threshold of roughly 0.5V. Thus, the output voltage at $n_2$ jumps up, if the sinusoidal input voltage an $n_1$ exceeds a certain value and goes down slowly due to the discharging of the capacitor. The refinement algorithm was used with a refinement
rate $\alpha=2.5$, which is large enough to increase the grid size fast enough to stop the algorithms after a few step, while a larger $\alpha$ would result in a larger grid, which makes the final steps of the algorithm to time consuming.
The stopping tolerance was chosen as $\varepsilon=0.02$ so that the algorithm stops when no visual improvements can 
be seen anymore. It can be observed in Fig.~\ref{rect} that the grid is in particular refined at the sharp transient of the output signal near 25$\mu$s, leading to an improved representation of that signal. 

Finally the spline expansion is coarsened to reduce the size of the expansion, using a threshold of $0.006$ over several 
iterations. Obviously the coarsening makes the approximation much more efficient. This is due to the fact the wavelet based
local error estimator yields only a guess of the refinement location. Although, this is a relative good guess, it is still a guess,
and the solution can be essentially improved by throwing out  unnecessary spline knots after the simulation.

\subsection{Multirate simulation}

Although spline wavelet based circuit simulation yields reliable results at reasonable cost
for initial value or periodic steady state problems, its most
interesting application is in multirate envelope simulation. Details about multirate circuit simulation can be found
in \cite{Bra2001,Hou04,Pulch08b,Pulch08a,BiBra12b,BiBra14b}. 
Here we follow the presentation in \cite{BiBra12b,BiBra14b}.

In multirate problems, the solution can be described as a fast oscillating carrier signal modulated by a slow changing
envelope signal, as it occurs usually in Radio Frequency (RF) problems. To separate the different scales 
the ordinary DAE's (\ref{mna}) is replaced by a system of partial
differential algebraic equations (PDAE's) of the form 
\begin{equation}\label{multirate}
\tfrac{\partial}{\partial \tau} q\big(\hat{x}(\tau,t)\big)
+\omega(\tau)\,\tfrac{\partial}{\partial t} q\big(\hat{x}(\tau,t)\big)
+g\big(\hat{x}(\tau,t)\big)=\hat{s}\big(\tau,t\big)
\end{equation}
under the conditions
$$
\hat{x}(0,t)=X_0(t),\quad  \hat{x}(\tau,t)=\hat{x}(\tau,t+P).
$$
The solution of (\ref{mna}) is obtained as the solution of (\ref{multirate})
along characteristic curves, i.e., $x(t)=\hat{x}\big(t,\Omega(t)\big)$. 
For typical multirate signals the bivariate envelope solution $\hat{x}(\tau,t)$ is smooth such that is
can be approximated
much more efficient than the original solution $x(t)$, since the limits of Nyquist's sampling theorem are avoided.

For the numerical solution we perform a semi-discretization with respect to $\tau$ (Rothe method), 
which is done by a multi-step method (namely Gear's BDF technique, see e.g.~\cite{hairer1993solving,hairer2010solving}). 
By this approach we find 
for each time step $\tau_k$, $k\in\NN$, an
approximation $X_k(t)$  of $\hat{x}(\tau_k,t)$, $\tau_{k+1}>\tau_k>\tau_0=0$,  as solution of an
ordinary DAE
\begin{equation}\label{envelope}
\omega(\tau_k)\tfrac{d}{dt}q\big(x(t)) + g_k(x(t),t) = 0,\qquad x(t) = x(t+P).
\end{equation}
Here, $g_k(x,t)$ 
is determined by the used multi-step method and depends therefore on solutions $X_\ell$, $\ell<k$,
at previous time steps (see \cite{BiBra12b,BiBra14b} for details).

The periodic problem (\ref{envelope}) is of the same structure as the original circuit equations (\ref{mna}) and
thus solved by the method described in \S\ref{pss}.  For typical
multirate problems the solution $\hat{x}(\tau,t)$ is very smooth in $\tau$.
Thus, the solution $X_{k-1}$ of the previous step can be expected to be a good approximation
for $X_k$ and therefore a good initial guess for the first Newton iteration (Step 1.\  in  Algo.~\ref{refine}). 
But due to the grid refinements the size of the spline grid would increase for each time step $\tau_k$.
Therefore a wavelet based grid  coarsening (\S\ref{coarse}) is applied to $X_{k-1}$, to get rid of unnecessary 
knots. Summing up, the \emph{refinements} will provide the required accuracy, 
while the grid \emph{coarsening} is responsible for the
sparsity of the grid, which ensures the efficiency of the simulation.

One advantage of our method is that information from previous steps is used for grid generation and Newton's 
initial guess. This is a crucial difference to an earlier approach in \cite{BaKnPu09}, where the \emph{uniform} 
wavelet decomposition of the solution of an initial value problem is used to generate an adaptive grid for a finite 
difference method on the periodic problem. Due to an excellent initial guess, Newton's
method  stops after a few iteration steps, often needing only one refinement with a moderate refinement rate
$\alpha < 2$. This can be achieved with a properly tuned step size control for the BDF method, and a threshold
for the coarsening, which is an order of magnitude below the required accuracy.
In our algorithm we use a time step
control, which uses a priori and a posteriori error estimates and simultaneously controls the number of Newton steps. The optimization of
this time step control is the object of current research.  

\paragraph*{Phase Locked Loop (PLL)}
We have used the described method to simulate a Phase Locked Loop (PLL) containing 205 MOSFET Transistors, and 
$N=112$ unknowns. This is a relatively complex circuit consisting of a voltage controlled oscillator (VCO),
frequency dividers, a phase frequency detector (PFD), and a loop filter.
For details we refer to \cite{BiBra12b,BiBra14b}. The input signal is a frequency modulated sinusoidal signal with center frequency 25kHz.
The baseband signal is also sinusoidal with frequency 10Hz and frequency deviation 100Hz. Central components as frequency dividers and the PFD are digital circuitry so that many internal signals exhibit sharp transients. Thus, adaptive grid 
generation shall lead to improved performance. We have chosen $P=40\mu s$ which corresponds to the center frequency.
The factor $\omega(\tau_k)$ is chosen by a method described in \cite{BiBra12b} to get a smooth solution. 

\begin{figure}[htb]
\parbox[t]{0.48\linewidth}{
\includegraphics[width=\linewidth]{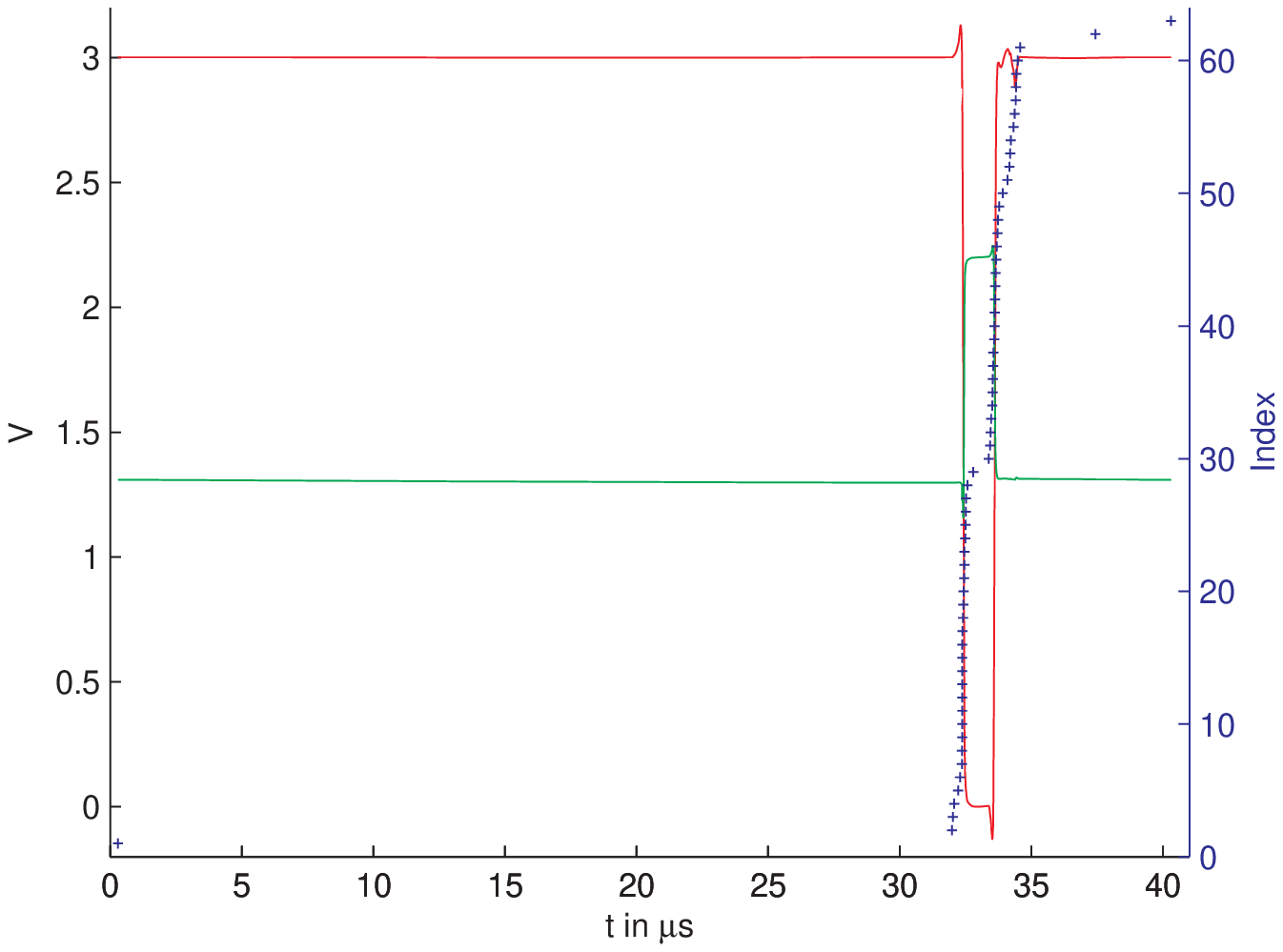}
\caption{\label{pll20} Two simulated signals of the PLL and the corresponding spline knots $t_k$, where each $+$ stands for a pair
$(t_k,k)$.}
}\hfill
\parbox[t]{0.48\linewidth}{
\includegraphics[width=\linewidth]{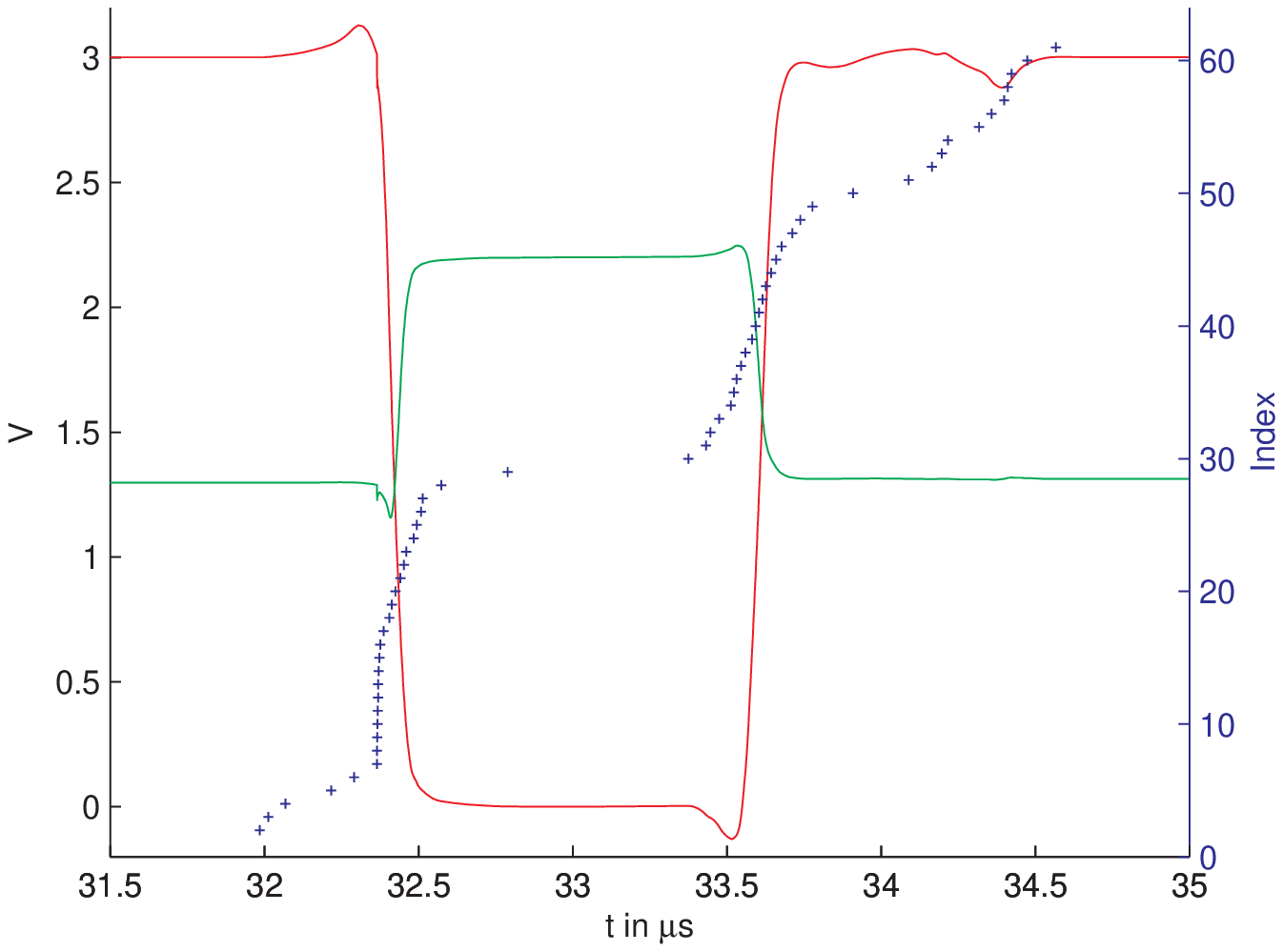}
\caption{\label{pll20det} Detail of Fig.~\ref{pll20}.}
}
\end{figure}

Fig.~\ref{pll20} shows two components of the approximated solution $X_{20}$ at $\tau_{20} = 72\mu s$ 
together with the corresponding spline grids. For each spline knot we have plotted the pair $(t_k,k)$, which gives a better
picture of the grid in particular at locations where the grid is locally dense. Together with the detail plot in Fig.~\ref{pll20det}
the example shows an excellent adaptation to the signal shapes, with high knot density at locations of sharp signal transients.

\begin{figure}[htb]
\hspace*{-1.8em}
\includegraphics[height=31.5ex]{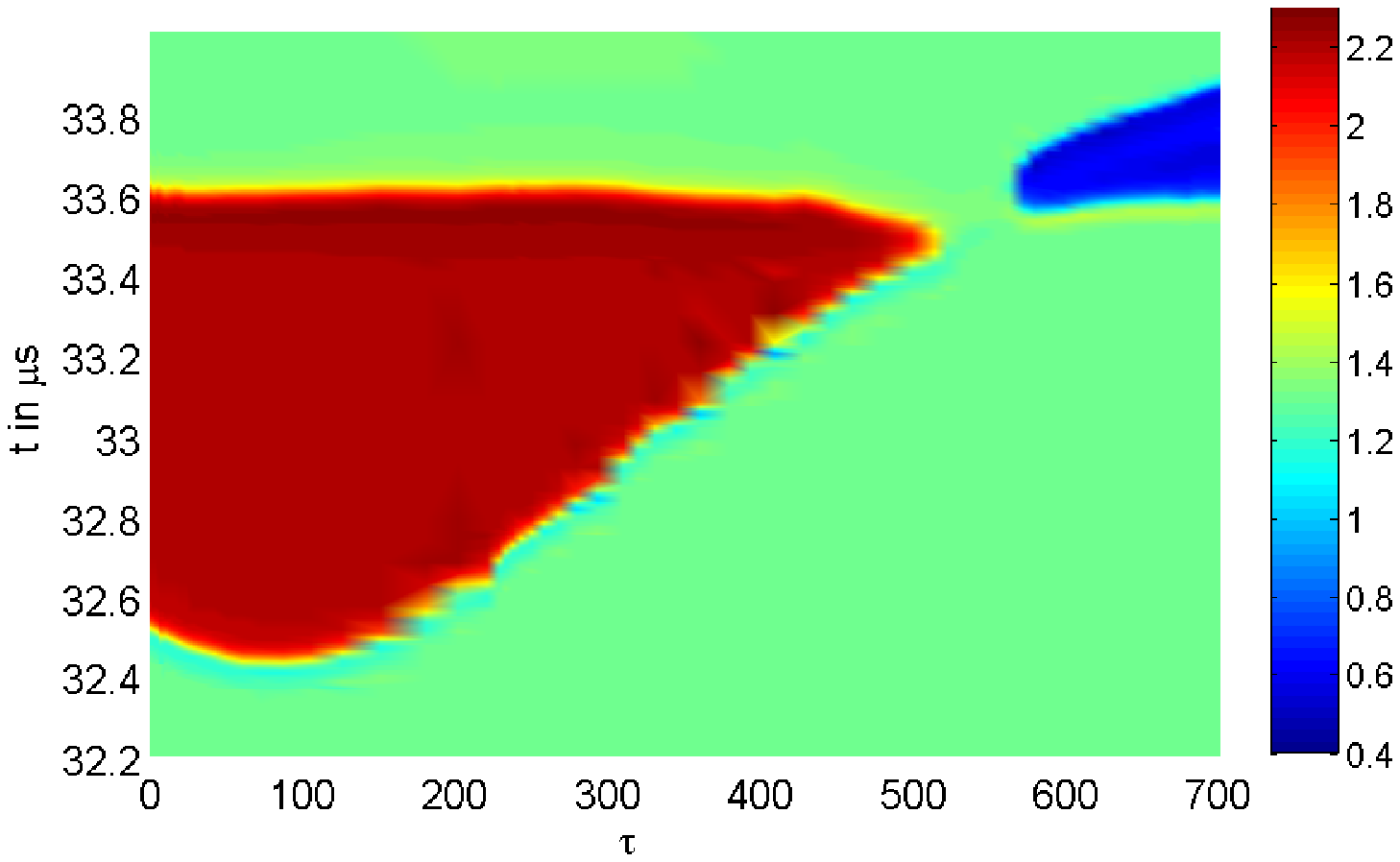}
\hspace*{-2.5em}
\includegraphics[height=31.5ex]{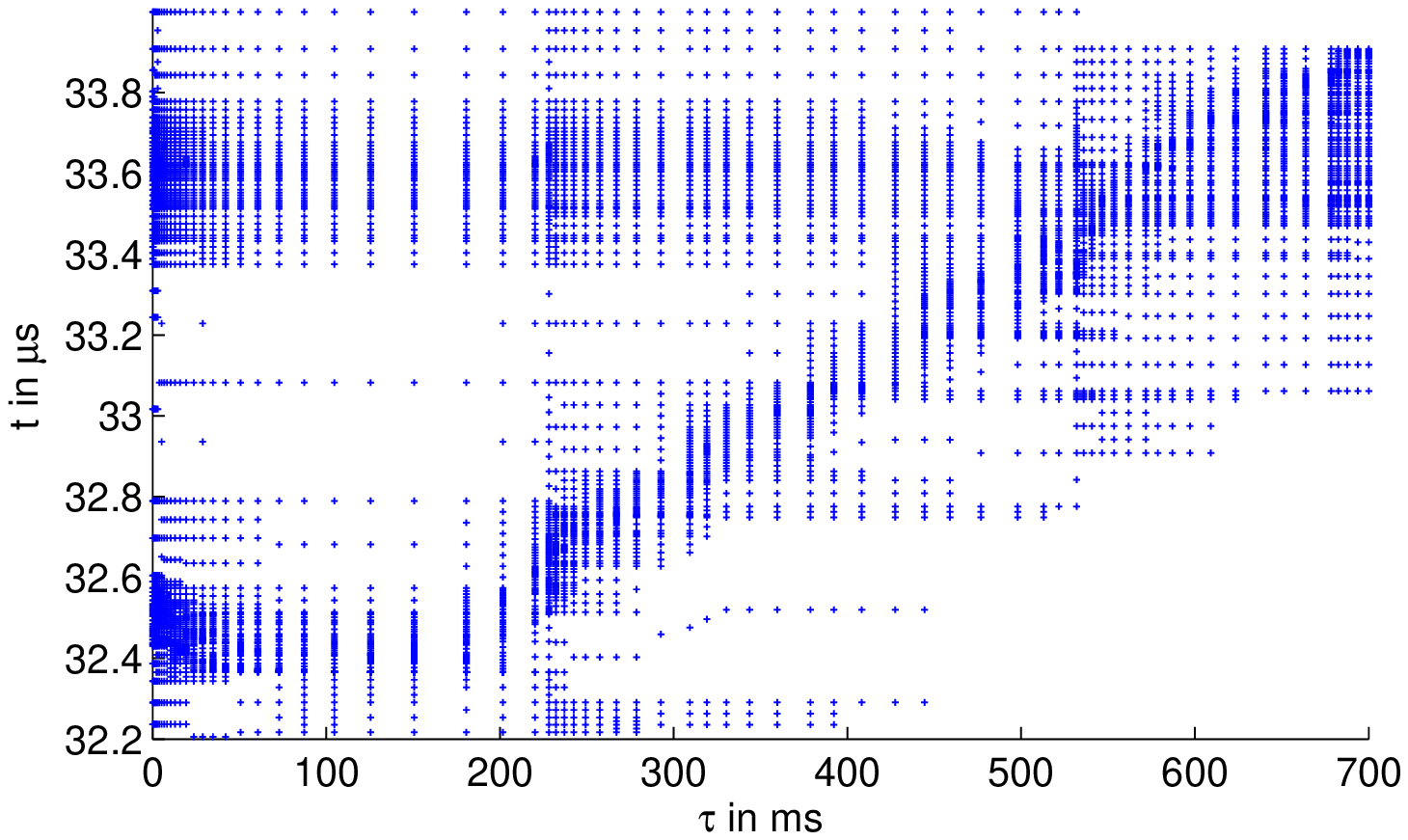}

\caption{\label{pll_pump} Detail of the solution of the multirate simulation  of a PLL (density plot)
together with  the discretization grid.}
\end{figure}

The grid development over several time steps $\tau_k$ is illustrated in Figure~\ref{pll_pump}. 
One can see that the grid is refined at the location of sharp signal transients, while the coarsening successfully removes knots
which are no longer needed in later time steps.

\paragraph*{Collpitts Oscillator}
\begin{figure}[htb]
\centerline{\includegraphics[width=0.4\linewidth]{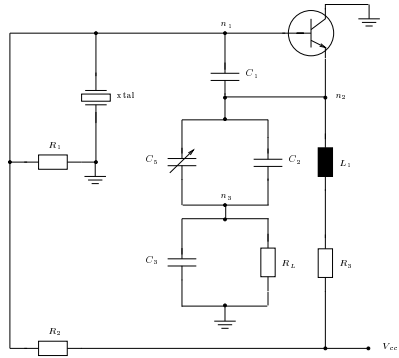}}
\caption{3MHz Colpitts quartz oscillator.\label{colp}}
\end{figure}

Another example is the start up phase of a  3MHz Colpitts quartz crystal oscillator (see Fig.~\ref{colp}). 
With $N=8$ this is a smaller
circuit. More details can be found in \cite{BiBra12a}, where also the problem of numerical damping is 
treated. For this problem we have chosen $P=1$, so that $\omega(\tau_k)$ determined by the method 
from \cite{BiBra12b} will be close to $3\cdot 10^6$, which corresponds to the oscillator frequency.

\begin{figure}
 \begin{center}
 \includegraphics[width=0.49\columnwidth]{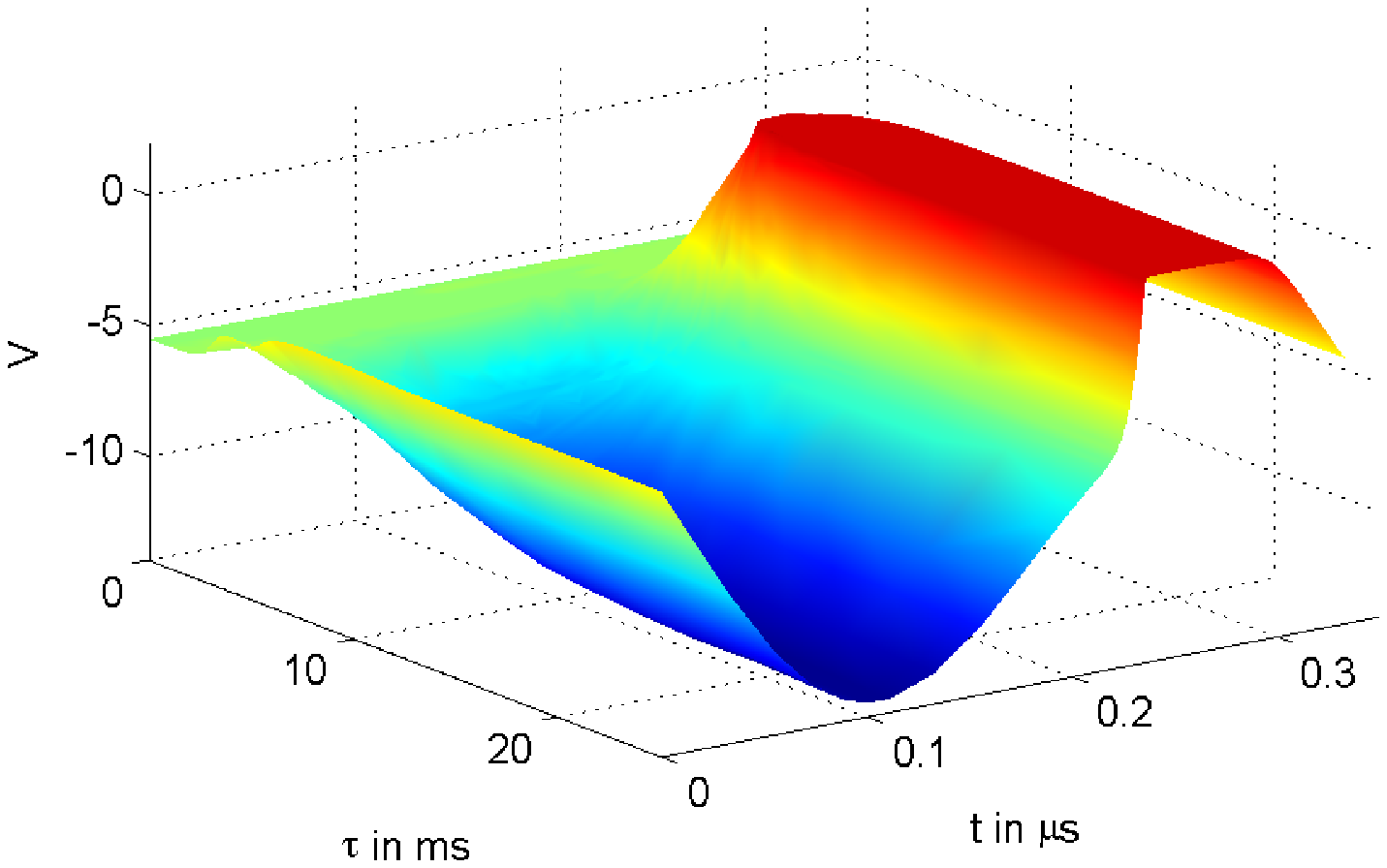}
 \includegraphics[width=0.49\columnwidth]{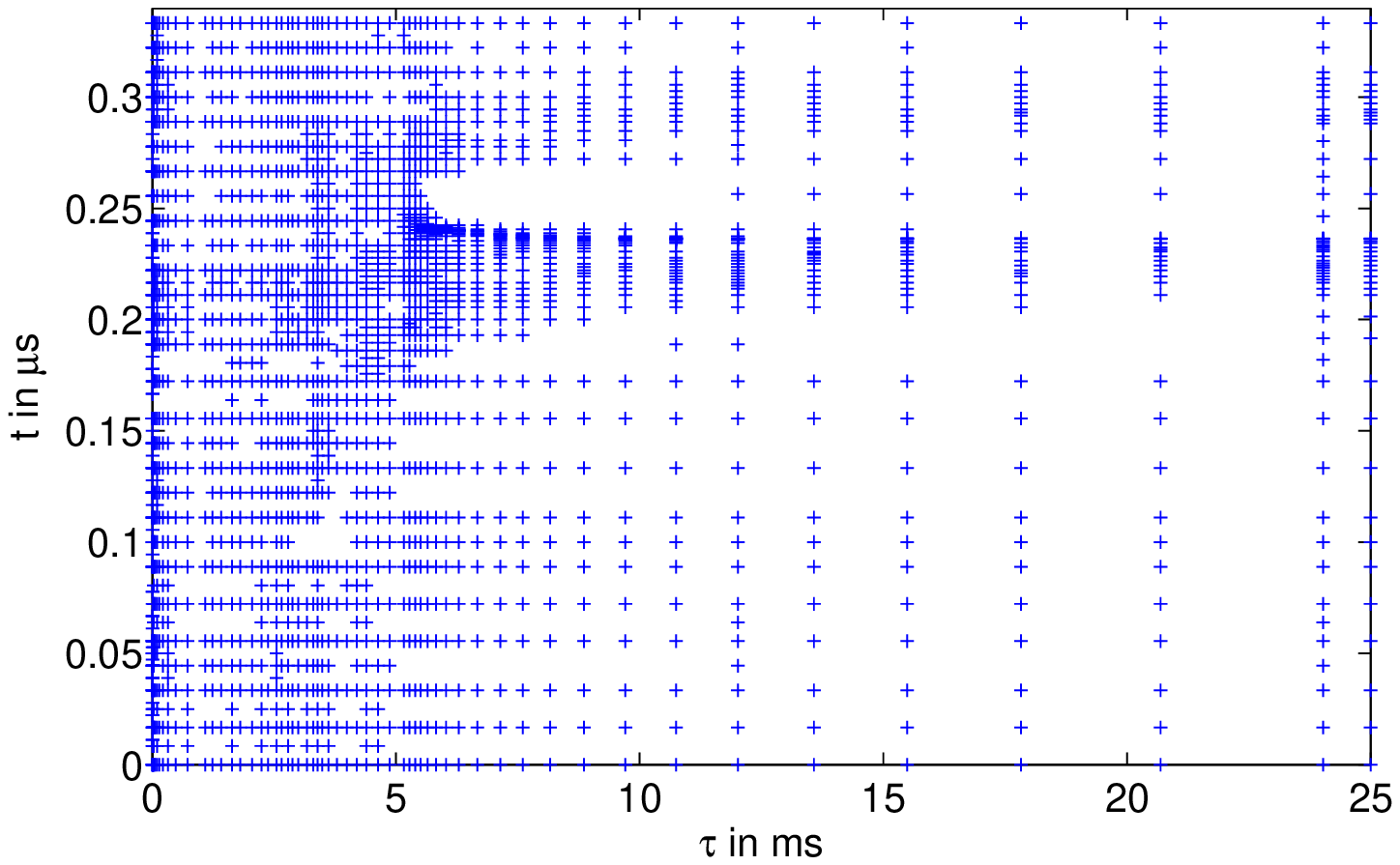}
 \caption{\label{oscillator} Multi-rate solution for 3MHz Colpitts oscillator and the corresponding adaptive grid.}
 \end{center}
\end{figure}

The transistor used as feedback amplifier introduces some nonlinear effects into the output signal
(see Fig.~\ref{oscillator}), which result in a sharp edge near $t=0.24\mu s$, which becomes apparent at the end of the 
start up phase  at $\tau\approx 5ms$. Obviously the grid refinement can handle this emerging edge, and adapts even
to a change of location for increasing $\tau$.  

\bibliographystyle{siam}
\bibliography{../bib/wavelet,../bib/other,../bib/biblio}

\end{document}